\documentclass{amsart}          

\usepackage{amsmath,amsthm}     
\usepackage{amssymb}            
\usepackage{euscript}           
\usepackage{enumerate,calc}
\usepackage[matrix,arrow,curve,frame]{xy}    

\xymatrixcolsep{1.9pc}                          
\xymatrixrowsep{1.9pc}
\newdir{ >}{{}*!/-5pt/\dir{>}}                  

\def\longfib{\DOTSB\relbar\joinrel\twoheadrightarrow}

\newtheorem{thm}[subsection]{Theorem}
\newtheorem{defn}[subsection]{Definition}
\newtheorem{prop}[subsection]{Proposition}

\newtheorem{cor}[subsection]{Corollary}
\newtheorem{lemma}[subsection]{Lemma}

\theoremstyle{definition}  
\newtheorem{example}[subsection]{Example}

\newtheorem{remark}[subsection]{Remark}

  {\end{list}}

\newcommand{\dfn}{\textbf} 

\newcommand{\mdfn}[1]{\dfn{\mathversion{bold}#1}} 

\newcommand{\Smash}             {\wedge}

\newcommand{\tens}              {\otimes}               
\newcommand{\iso}               {\cong}  

\newcommand{\cat}{\EuScript}    
\newcommand{\cA}{{\cat A}}      
\newcommand{\cC}{{\cat C}}
\newcommand{\cD}{{\cat D}}
\newcommand{\cE}{{\cat E}}
\newcommand{\cF}{{\cat F}}

\newcommand{\cI}{{\cat I}}

\newcommand{\cM}{{\cat M}}
\newcommand{\cN}{{\cat N}}

\newcommand{\cV}{{\cat V}}

\newcommand{\Set}{{\cat Set}}

\newcommand{\sSet}{s{\cat Set}}

\newcommand{\Ho}{{\cat Ho}}
\newcommand{\ho}{\text{Ho}\,}

\newcommand{\field}[1]  {\mathbb #1} 

\newcommand{\F}         {\field F}

\newcommand{\Z}         {\field Z}

\renewcommand{\S}       {\field S}

\DeclareMathOperator*{\colim}{colim}

\DeclareMathOperator{\Hom}{Hom}

\DeclareMathOperator{\Map}{Map}

\newcommand{\ra}{\rightarrow}                   
\newcommand{\lra}{\longrightarrow}              
\newcommand{\la}{\leftarrow}                    
\newcommand{\lla}{\longleftarrow}               
\newcommand{\llra}[1]{\stackrel{#1}{\lra}}      
\newcommand{\llla}[1]{\stackrel{#1}{\lla}}      

\newcommand{\we}{\llra{\sim}}                   
\newcommand{\bwe}{\llla{\sim}}
\newcommand{\cof}{\rightarrowtail}              
\newcommand{\trfib}{\stackrel{\sim}{\longfib}}
\newcommand{\trcof}{\stackrel{\sim}{\cof}}

\newcommand{\inc}{\hookrightarrow}              
\newcommand{\dbra}{\rightrightarrows}           

\newcommand{\blank}{-}                          
\newcommand{\Id}{Id}                            
\newcommand{\und}{\underline}
\newcommand{\mM}{\underline{\cM}}

\newcommand{\sing}{Sing}

\newcommand{\ovcat}{\downarrow}

\newcommand{\adjoint}{\rightleftarrows}

\newcommand{\he}{\simeq}

\numberwithin{equation}{section}

\newcommand{\cU}{{\cat U}}
\newcommand{\U}{\cU}
\newcommand{\Ua}{\cU_{ad}}
\newcommand{\Spe}{Sp}
\DeclareMathOperator{\hEnd}{hEnd}
\newcommand{\RingSpectra}{\cat RingSpectra}
\newcommand{\mD}{\underline{\cD}}
\newcommand{\mN}{\underline{\cN}}
\DeclareMathOperator{\Sing}{Sing}
\DeclareMathOperator{\Real}{Re}
\DeclareMathOperator{\Func}{Func}
\newcommand{\unit}{\mathbf{1}}
\DeclareMathOperator{\Ev}{Ev}
\DeclareMathOperator{\hdga}{hEnd_{ad}}
\newcommand{\Ab}{\cat Ab}
\newcommand{\Ad}{Ad}

\newcommand{\ME}{ME}

\newcommand{\Spectra}{\Spe^\Sigma}
\newcommand{\cofib}{\cof}
\renewcommand{\F}{F}
\newcommand{\Ring}{\cat{R}ing}

\begin{document}

\title{Spectral enrichments of model categories}

\author{Daniel Dugger}
\address{Department of Mathematics\\ University of Oregon\\ Eugene, OR
97403} 

\email{ddugger@math.uoregon.edu}

\begin{abstract}
We prove that every stable, combinatorial model category can be
enriched in a natural way over symmetric spectra.  
As a consequence of the
general theory, every object in such a model category has an
associated homotopy endomorphism ring spectrum.  Basic properties
of these invariants are established.
\end{abstract}

\maketitle

\tableofcontents

\section{Introduction}

If $X$ and $Y$ are two objects in a model category $\cM$, it is
well-known that there is an associated `homotopy function complex'
$\Map(X,Y)$ (cf. \cite[Chap. 17]{H} or \cite[Sec. 5.4]{Ho2}).  This is
a simplicial set, well-defined up to weak equivalence, and it is an
invariant of the homotopy types of $X$ and $Y$.  Following \cite{DK}
one can actually construct these function complexes so that they come
with composition maps $\Map(Y,Z) \times \Map(X,Y) \ra \Map(X,Z)$, thereby
giving an enrichment of $\cM$ over simplicial sets.  This enrichment
is an invariant (in an appropriate sense) of the model category $\cM$.

This paper concerns analagous results for stable model categories,
with the role of simplicial sets being replaced by symmetric spectra
\cite[Th. 3.4.4]{HSS}.  We show that if $\cM$ is a stable,
combinatorial model category then any two objects can be assigned a
symmetric spectrum function complex.  More importantly, one can give
composition maps leading to an enrichment of $\cM$ over the symmetric
monoidal category of symmetric spectra.  One application is that any
object $X\in \cM$ has an associated `homotopy endomorphism ring
spectrum' $\hEnd(X)$ (where by ring spectrum we mean essentially what
used to be called an $A_\infty$-ring spectrum).  These ring spectra,
as well as the overall enrichment by symmetric spectra, are homotopy
invariants of the model category $\cM$.

\medskip

\subsection{An application}
Before describing the results in more detail, here is the motivation
for this paper.  If $R$ is a differential graded algebra, there is a
stable model category structure on (differential graded) $R$-modules
where the weak equivalences are quasi-isomorphisms and the fibrations
are surjections.  Given two dgas $R$ and $S$, when are the model
categories of $R$- and $S$-modules Quillen equivalent?  A complete
answer to this question is given in \cite{DS}.  The problem is subtle:
even though the categories of $R$- and $S$-modules are additive,
examples show that it's possible for them to be Quillen equivalent only
through a zig-zag involving non-additive model categories.  To deal
with this, the arguments in \cite{DS} depend on using homotopy
endomorphism ring spectra as invariants of stable model categories.
The present paper develops some of the tools necessary for those
arguments.

\subsection{Statement of results}
A category is {\it locally presentable\/} if it is cocomplete and all
objects are small in a certain sense; see \cite{AR}.  A model category
is called \dfn{combinatorial} if it is cofibrantly-generated and the
underlying category is locally presentable.  This class was introduced
by Jeff Smith, and the examples are ubiquitous (the class even
includes model categories made from topological spaces, if one uses
$\Delta$-generated spaces).  Background information on combinatorial
model categories can be found in \cite{D2}.

A model category is called \dfn{stable} if the initial and terminal
objects coincide (that is, it is a {\it pointed\/} category) and if
the induced suspension functor is invertible on the homotopy category.

Our results concern enrichments of stable, combinatorial model
categories.  Unfortunately we do not know how to give a {\it
canonical\/} spectral enrichment for our model categories, however.
Instead there are many such enrichments, involving choices, but the
choices yield enrichments which are homotopy equivalent in a certain
sense.  The machinery needed to handle this is developed in
Section~\ref{se:modelenrich}.  There we define a \dfn{model
enrichment} of one model category by another, and give a notion of two
model enrichments being \dfn{quasi-equivalent}.  A crude version of
our main theorem can be stated as follows:

\begin{thm}
\label{th:spenrich}
Every stable, combinatorial model category has a canonical
quasi-equivalence class of model
enrichments by $\Spectra$.
\end{thm}

Here $\Spe^\Sigma$ denotes the model category of symmetric spectra from
\cite{HSS}, with its symmetric monoidal smash product.  `Canonical'
means the enrichment has good functoriality properties with respect to
Quillen pairs and Quillen equivalences.  More precise statements are
given in Section~\ref{se:main}.  We will show that the
canonical enrichment by $\Spectra$ is preserved, up to
quasi-equivalence, when you prolong or restrict across a Quillen
equivalence.  It follows that the enrichment contains only
`homotopy information' about the model category; so it can be used to
decide whether or not two model categories are Quillen equivalent.

One simple consequence of the above theorem is the following:

\begin{cor}
\label{co:hoenrich}
If $\cM$ is a stable, combinatorial model category then $\Ho(\cM)$ is
naturally enriched over $\Ho(\Spectra)$.  
\end{cor}

The above corollary is actually rather weak, and not representative of
all that the theorem has to offer.  For instance, the corollary
implies that every object of such a model category has an
endomorphism ring object in $\Ho(\Spectra)$---that is, a spectrum
$R$ together with a pairing $R\Smash R \ra R$ which is associative and
unital up to homotopy.  The theorem, on the other hand, actually gives
the following:

\begin{cor}
Every object $X$ of a stable, combinatorial model category has a
naturally associated $A_\infty$-ring spectrum $\hEnd(X)$---called the
\dfn{homotopy endomorphism spectrum} of $X$---well-defined in the
homotopy category of $A_\infty$-ring spectra.  If $X\he Y$ then
$\hEnd(X)\he \hEnd(Y)$.  
\end{cor}

The main results concerning these endomorphism spectra are as follows.
The first shows that they are homotopical invariants of the model
category $\cM$:

\begin{thm}
\label{th:invariance}
Let $\cM$ and $\cN$ be stable, combinatorial model categories.
Suppose they are Quillen equivalent, through a zig-zag where the
intermediate steps are possibly not combinatorial or pointed.  Let
$X\in \cM$, and let $Y\in \ho(\cN)$ be the image of $X$ under the
derived functors of the Quillen equivalence.  Then $\hEnd(X)$ and
$\hEnd(Y)$ are weakly equivalent ring spectra.
\end{thm}

A model category $\cM$ is called \dfn{spectral} if it is enriched,
tensored, and cotensored over symmetric spectra in a homotopically
well-behaved manner ($\cM$ is also called an $\Spe^\Sigma$-model
category).  See Section~\ref{se:spectral} for a more detailed
definition.  The following result says that in spectral model
categories homotopy endomorphism spectra can be computed in the
expected way, using the spectrum hom-object
$\und{\cM}_{\Spe^\Sigma}(\blank,\blank)$:

\begin{prop}
\label{pr:spectralcase}
Let $\cM$ be a stable, combinatorial model category which is also
spectral.  Let $X$ be a cofibrant-fibrant object of $\cM$.
Then $\hEnd(X)$ and $\und{\cM}_{\Spe^\Sigma}(\blank,\blank)$
 are weakly equivalent ring spectra.
\end{prop}

Enhanced results are proven in the case where $\cM$ is also an {\it
additive\/} model category (see Section~\ref{se:additive} for the
definition).  In this context one obtains an enrichment over the
monoidal model category $\Spe^\Sigma(s\Ab)$ of symmetric spectra based
on simplicial abelian groups.  The paper \cite{S} shows this category
is monoidally equivalent to the model category of unbounded chain
complexes of abelian groups, which perhaps is easier to think about.
Any object $X\in \cM$ therefore has an {\it
additive homotopy endomorphism ring object\/} in $\Spe^\Sigma(s\Ab)$
(or equivalently, a ``homotopy endomorphism dga'').  These
endomorphism dgas are invariant under Quillen equivalences between
additive model categories, but not general Quillen equivalences.
Details are in Section~\ref{se:additive}.
  
\vspace{0.1in}

\subsection{The construction}
In \cite{DK} Dwyer and Kan constructed model enrichments over $\sSet$
via their hammock localization.  This is a very elegant construction,
in particular not involving any {\it choices\/}.  Unfortunately we
have not been clever enough to find a similar construction for
enrichments by symmetric spectra.  The methods of the present paper
are more of a hack job: they get us the tools we need at a relatively
cheap cost, but they are not so elegant.

The idea is to make use of the `universal' constructions from
\cite{D1,D2}, together with the general stabilization machinery
provided by \cite{Ho1}.  By \cite{D2} every combinatorial model
category is Quillen equivalent to a localization of diagrams of
simplicial sets.  Using the simplicial structure on this diagram
category, we can apply the symmetric spectra construction of
\cite{Ho1}.  This gives a new model category, Quillen equivalent to
what we started with, where one has actual symmetric spectra function
complexes built into the category.

In more detail, given a pointed, combinatorial model category
$\cM$ one can choose a Quillen equivalence $\U_+\cC/S \we \cM$ by
modifying the main result of \cite{D2}.  Here $\U_+\cC$ is the
universal {\it pointed\/} model category built from $\cC$, developed
in Section~\ref{se:univ}; $S$ is a set of maps in $\U_+\cC$, and
$\U_+\cC/S$ denotes the Bousfield localization \cite[Sec. 3.3]{H}.

The category $\U_+\cC/S$ is a nice simplicial model category, and we
can form symmetric spectra over it using the results of \cite{Ho1}.
This gives us a new model category $\Spe^\Sigma(\U_+\cC/S)$, which is
enriched over $\Spectra$.  If $\cM$ was stable to begin with then we
have a zig-zag of Quillen equivalences
\[ \cM \bwe \U_+\cC/S \we \Spe^\Sigma(\U_+\cC/S) 
\]  
and can transport the enrichment of the right-most model category onto
$\cM$.  Finally, theorems from \cite{D1} allow us to check that the
resulting enrichment of $\cM$ doesn't depend (up to quasi-equivalence)
on our chosen Quillen equivalence $\U_+\cC/S\ra \cM$.

\vspace{0.1in}

By now the main shortcoming of this paper should be obvious: all the
results are proven only for {\it combinatorial\/} model categories.
This is an extremely large class, but it is very plausible that the
results about spectral enrichments hold in complete generality.
Unfortunately we have not been able to find proofs in this
setting, so it remains a worthwhile challenge.

\vspace{0.2in}

\subsection{Organization of the paper}
Sections 2 and 3 contain the basic definitions of enrichments, model
enrichments, and the corresponding notions of equivalence.  Section 4
deals with the universal pointed model categories $\U_+\cC$, and
establishes their basic properties.  The main part of the paper is
Section 5, which gives the results on spectral enrichments and
homotopy endomorphism spectra.  Section 6 returns to the
proof of Proposition~\ref{pr:equiv=direct}: this is a foundational
result showing that quasi-equivalent enrichments have the properties
one hopes for.  Finally, in Section 7 we present expanded results for
additive model categories.  This entails developing
`universal additive model categories', a topic which may be of
independent interest.

We also give two appendices.  Appendix A contains several basic
results about model categories which are enriched, tensored, and
cotensored over a monoidal model category (the main examples for us
are {\it simplicial\/} and {\it spectral\/} model categories).  The
reader is encouraged to familiarize himself with this section before
tackling the rest of the paper.  Appendix B gives a general result
about commuting localization and stabilization.
   
\subsection{Terminology}
We assume a familiarity with model categories and localization theory,
for which \cite{H} is a good reference.  Several conventions from
\cite{D1} are often used, so we'll now briefly recall these.  A {\it
Quillen\/} map $L\colon \cM \ra \cN$ is another name for a Quillen
pair $L\colon \cM \adjoint \cN \colon R$.  If $L_1$ and $L_2$ are two
such Quillen maps, a {\it Quillen homotopy\/} $L_1 \ra L_2$ is a
natural transformation between the left adjoints which is a weak
equivalence on cofibrant objects.  If $\cM$ is a model category and
$S$ is a set of maps in $\cM$, then $\cM/S$ denotes the left Bousfield
localization \cite[Sec. 3.3]{H}.  

\vspace{0.1in}

\subsection{Acknowledgments}
Readers should note that the present paper owes a great
debt to both \cite{Ho1} and \cite{SS2}. 

\vspace{0.2in}


\section{Enrichments in category theory}

In this section we review the notion of a category being
\dfn{enriched} over a symmetric monoidal category.  Our situation is
slightly more general than what usually occurs in the
literature.  There is a notion of \dfn{equivalence} which encodes when
two enrichments carry the same information.

\medskip

\subsection{Basic definitions}
Let $\cC$ be a category, and let $(\cD,\tens,S)$ be a symmetric
monoidal category (where $S$ is the unit).  An \dfn{enrichment} of $\cC$ by
$\cD$ is a functor $\tau \colon \cC^{op} \times \cC \ra \cD$ together
with
\begin{enumerate}[(i)]
\item For every $a,b,c\in \cC$ a `composition map' $\tau(b,c)\tens
\tau(a,b) \ra \tau(a,c)$, natural in $a$ and $c$;
\item a collection
of maps $S \ra \tau(c,c)$ for every $c\in \cC$.
\end{enumerate}
This data is required to satisfy the associativity and unital rules
for composition, which are so standard that we will not write them
down.  We also require that for any map $f\colon a\ra b$ in $\cC$,
the square
\[ \xymatrix{ S \ar[r] \ar[d] & \tau(a,a) \ar[d]^{{f}} \\
              \tau(b,b) \ar[r]^{f} & \tau(a,b)
}
\]
commutes.  

Note that if $\cC=\{*\}$ is the trivial category and $\Ab$ is
the category of abelian groups, then an enrichment of $\cC$ by $\Ab$
is just another name for an associative and unital ring.

If $\tau$ and $\tau'$ are two enrichments of $\cC$ by $\cD$, a map
$\tau \ra \tau'$ is a natural transformation $\tau(a,b) \ra
\tau'(a,b)$ compatible with the unit and composition maps.

\begin{remark}
The above definition differs somewhat from related things in the
literature.  According to \cite[Sec. 6.2]{B}, a \mdfn{$\cD$-category}
is a collection of objects $\cI$ together with a Hom-object
$\cI(i,j)\in \cD$ for every $i,j\in \cI$, etc.  This corresponds to
our above definition in the case where $\cC$ has only identity maps.

If $\cC$ is a category (i.e., a $\Set$-category), one can define a
$\cD$-category $\S_\cC$ with the same object set as $\cD$ and 
$\S_\cC(a,b)=\coprod_{\cC(a,b)} S$.  To
give an enrichment of $\cC$ by $\cD$ in the sense we defined above is
the same as giving a $\cD$-category with the same objects as $\cC$,
together with a $\cD$-functor from $\S_\cC$ to this $\cD$-category.
\end{remark}

\begin{example}
If $\cM$ is a simplicial model category, the assignment $X,Y \mapsto
\Map(X,Y)$ is an enrichment of $\cM$ by $\sSet$.  If $\cM$ is a
general model category, the hammock localization assignment $X,Y
\mapsto L^H\cM(X,Y)$ from 
\cite[3.1]{DK} is also an enrichment of $\cM$ by $\sSet$.
\end{example}

\subsection{Bimodules}
\label{se:bimodules}
Let $\sigma$ and $\tau$ be two enrichments of $\cC$ by $\cD$.  By a
$\sigma-\tau$ \dfn{bimodule} we mean a collection of objects
$M(a,b)\in \cD$ for every $a,b\in \cC$, together with `multiplication
maps'
\[ \sigma(b,c) \tens M(a,b) \ra M(a,c) \qquad M(b,c)\tens \tau(a,b)
\ra M(a,c) 
\]
which are natural in $a$ and $c$.  We again assume associativity and
unital conditions which we will not write down, as well as the
property that for any $a,b,c,d\in \cC$ the two obvious maps
\[ \sigma(c,d) \tens M(b,c) \tens \tau(a,b) \dbra M(a,d)
\]
are equal.

Note that a bimodule has a natural structure of a bifunctor
$\cC^{op}\times \cC \ra \cD$.  For instance, if $f\colon a\ra b$ is a
map in $\cC$ then consider the composite $S \ra \sigma(a,a) \ra
\sigma(a,b)$.  We then have $S\tens M(a',a) \ra \sigma(a,b)\tens
M(a',a) \ra M(a',b)$, giving a map $M(a',a) \ra M(a',b)$ induced by
$f$.  Similar considerations give functoriality in the first variable.

\begin{remark}
For a more precise version of the definition of bimodule, see
Section~\ref{se:bimodules2}.  Earlier parts of
Section~\ref{se:leftover} also define the notions of left and right
$\sigma$-module, which we have for the moment skipped over.
\end{remark}

To understand the following definition, observe that two rings $R$ and
$S$ are isomorphic if and only if there is an $R-S$ bimodule $M$
together with a chosen element $m\in M$ such that the induced maps
$r \ra rm$ and $s \ra ms$ give isomorphisms of abelian groups 
$R \ra M \la S$.

\begin{defn}
\label{de:equiv}
Let $\sigma$ and $\tau$ be two enrichments of $\cC$ by $\cD$.  
\begin{enumerate}[(a)]
\item
By a \mdfn{pointed $\sigma-\tau$ bimodule} we mean a bimodule $M$
together with a collection of maps $S\ra M(c,c)$ for every $c\in \cC$,
such that for any map $a\ra b$ the square
\[ \xymatrix{ S \ar[r] \ar[d] & M(a,a) \ar[d] \\
              M(b,b) \ar[r] & M(a,b)
}
\]
commutes.
\item
We say that $\sigma$ and $\tau$ are \dfn{equivalent} if there is a
pointed $\sigma-\tau$ bimodule $M\colon \cC^{op}\times \cC \ra \cD$
for which the composites
\begin{align*}
\sigma(a,b)\tens S \ra \sigma(a,b) &\tens M(a,a) \ra M(a,b),  
\qquad\text{and} \\
 S\tens \tau(a,b) \ra M(b,b) &\tens \tau(a,b) \ra M(a,b)
\end{align*}
are isomorphisms, for every $a,b\in \cC$.  
\end{enumerate}
\end{defn}

\begin{remark}
A $\sigma-\tau$ bimodule is, by restriction, an $\S_\cC-\S_\cC$
bimodule.  Note that $\S_\cC$ has an obvious structure of
$\S_\cC-\S_\cC$ bimodule.  The definition of {\it pointed\/}
$\sigma-\tau$ bimodule says that there is a map of $\S_\cC-\S_\cC$
bimodules $\S_\cC \ra M$.
\end{remark}

\begin{lemma}
\label{le:equiv=iso}
Assume that $\cD$ has pullbacks.
Two enrichments $\sigma$ and $\tau$ are equivalent if and only if
there is an isomorphism $\sigma \iso \tau$.
\end{lemma}

\begin{proof}
If there is an isomorphism $\sigma\iso \tau$, then we let $M=\tau$ and
regard it as a $\sigma-\tau$ bimodule.  This shows $\sigma$ and $\tau$
are equivalent.

If we instead assume that $\sigma$ and $\tau$ are equivalent via the
pointed bimodule $M$, define $\theta(a,b)$ to be the pullback
\[ \xymatrix{ \theta(a,b) \ar[r] \ar[d] & \tau(a,b) \ar[d] \\
         \sigma(a,b) \ar[r] & M(a,b).
}
\]
Here the lower horizontal map is the composite 
\[ \sigma(a,b) \iso \sigma(a,b) \tens S \ra \sigma(a,b)\tens M(a,a)
\ra M(a,b)
\]
and the right vertical map is defined similarly.  The universal
property of the pullback allows one to see that $\theta$ is naturally
an enrichment of $\cC$ by $\cD$, and that $\theta \ra \sigma$ and
$\theta \ra \tau$ are maps of enrichments.

Now, our assumption that $\sigma$ and $\tau$ are equivalent via $M$
includes the condition that the bottom and right maps in the above
pullback square are isomorphisms.  So all maps in the square are
isomorphisms, which means we have $\sigma \iso \theta \iso \tau$.   
\end{proof}

\begin{remark}
Since the notions of equivalence and isomorphism coincide, one might
wonder why we bother with the former.  The answer is in the next
section, where the homotopical analogs of these two notions slightly
diverge.
\end{remark}


\section{Enrichments for model categories}
\label{se:modelenrich}

We now give model category analogs for the material from the last
section.  There is the notion of \dfn{model enrichment}, together with
two notions of equivalence: these are called \dfn{quasi-equivalence}
and \dfn{direct equivalence}.  Direct equivalences have the property
of obviously preserving the `homotopical' information in an
enrichment; but quasi-equivalences are what seem to arise in practice.
Fortunately the two notions are closely connected---see
Proposition~\ref{pr:equiv=direct}.

The material in this section is a simple extension of techniques from
\cite{SS2}, which dealt with enrichments over symmetric spectra.

\medskip

\subsection{Model enrichments}
Let $\cM$ be a model category and let $\cV$ be a symmetric monoidal
model category \cite[Def. 4.2.6]{Ho2}.  A \dfn{model enrichment} of
$\cM$ by $\cV$ is an enrichment $\tau$ with the property that whenever
$a\ra a'$ is a weak equivalence between cofibrant objects, and $x\ra
x'$ is a weak equivalence between fibrant objects, then the induced
maps
\[ \tau(a',x) \ra \tau(a,x) \qquad\text{and}\qquad
   \tau(a,x) \ra \tau(a,x')
\]
are weak equivalences.  

A \dfn{quasi-equivalence} between two model enrichments $\sigma$ and
$\tau$ consists of a pointed $\sigma-\tau$ bimodule $M$ such that the
compositions
\begin{align*}
\sigma(a,b)\tens S \ra \sigma(a,b) &\tens M(a,a) \ra M(a,b)
\qquad\text{and} \\
 S\tens \tau(a,b) \ra M(b,b) &\tens \tau(a,b) \ra M(a,b)
\end{align*}
are weak equivalences whenever $a$ is cofibrant and $b$ is fibrant.
 
\begin{defn}
Let \mdfn{$\ME_0(\cM,\cV)$} be the collection of equivalence classes of
model enrichments, where the equivalence relation is the one generated
by quasi-equivalence. 
\end{defn}

\begin{example}
\label{ex:example}
Let $\cM$ be a simplicial model category, and let $\tau(X,Y)$ be the
simplicial mapping space between $X$ and $Y$.  This is a model
enrichment of $\cM$ by $\sSet$.  Let $QX \trfib X$ be a
cofibrant-replacement functor for $\cM$, and define
$\tau'(X,Y)=\tau(QX,QY)$.  This is another model enrichment of
$\cM$, but note that there are no obvious maps between $\tau$ and
$\tau'$.  There is an obvious quasi-equivalence, however: define
$M(X,Y)=\Map(QX,Y)$.  This is a $\tau-\tau'$ bimodule, and the maps
$QX \ra X$ give the distinguished maps $* \ra M(X,X)$.

This example illustrates that quasi-equivalences arise naturally, more
so than the notion of `direct equivalence' we define  next.
\end{example}

\subsection{Direct equivalences}
A map of model enrichments $\tau \ra \tau'$ is a \dfn{direct
equivalence} if $\tau(a,b) \ra \tau'(a,b)$ is a weak equivalence
whenever $a$ is cofibrant and $b$ is fibrant.  

To say something about
the relationship between quasi-equivalence and direct equivalence, we
need a slight enhancement of our definitions.  If $\cI$ is a full
subcategory of $\cM$, we can talk about \mdfn{model enrichments
defined over $\cI$}: meaning that $\tau(a,b)$ is defined only
for $a,b \in \cI$.  In the same way we can talk about ``direct
equivalences over $\cI$'', and so on.

Now we can give the following analog of Lemma~\ref{le:equiv=iso}.  
This is the most important result of this section.  

\begin{prop}
\label{pr:equiv=direct}
Let  $\cV$ be a combinatorial, symmetric monoidal model category
satisfying the monoid axiom \cite[Def. 3.3]{SS1}.  Assume also that
the unit $S\in \cV$ is cofibrant.  Let $\sigma$ and $\tau$ be model
enrichments of $\cM$ by $\cV$.  Let $\cI$ be a small, full subcategory
of $\cM$ consisting of cofibrant-fibrant objects.  If $\sigma$ and
$\tau$ are quasi-equivalent over $\cI$, then there is a zig-zag of
direct equivalences (over $\cI$) between $\sigma$ and $\tau$.
\end{prop}

The assumption about the smallness of $\cI$ is needed so that there is
a model structure on certain categories of modules and bimodules, a
key ingredient of the proof.

\begin{proof}[Sketch of proof]
The proof can be adapted directly from \cite[Lemma A.2.3]{SS2}, which
dealt with the case where $\cV$ is symmetric spectra and $\cI$ has
only identity maps.  Essentially the proof is a homotopy-theoretic
version of the pullback trick in Lemma~\ref{le:equiv=iso}.  

Let $M$ be a bimodule giving an equivalence between $\sigma$ and
$\tau$.  When the maps $\sigma(a,b) \ra M(a,b)$ are trivial
fibrations, the pullback trick immediately gives a zig-zag of direct
equivalences between $\sigma$ and $\tau$.  For the general case
one uses certain model structures on module categories to reduce to
the previous case.  A full discussion requires quite a bit of
machinery, so we postpone this until Section~\ref{se:leftover}.
\end{proof}

\begin{cor}
\label{co:endo=invariant}
Let $\sigma$ and $\tau$ be model enrichments of $\cM$ by $\cV$.  Let
$X$ be a cofibrant-fibrant object of $\cM$.  If $\sigma$ and $\tau$
are quasi-equivalent, then the $\cV$-monoids $\sigma(X,X)$ and $\tau(X,X)$
are weakly equivalent in $\cV$ (meaning there is a zig-zag between
them where all the intermediate objects are monoids in $\cV$, and all
the maps are both monoid maps and weak equivalences).
\end{cor}

\begin{proof}
This is an application of Proposition~\ref{pr:equiv=direct}, where
$\cI$ is the full subcategory of $\cM$ whose sole object is $X$.
\end{proof}

\begin{cor}
\label{co:endo=coffib}
Let $\sigma$ be a model enrichment of $\cM$.  Let $\cI$ be a small
category, and let $G_1,G_2\colon \cI \ra \cM$ be two functors whose
images lie in the cofibrant-fibrant objects.  Assume there is a
natural weak equivalence $G_1 \we G_2$.  Then the enrichments on $\cI$
given by $\sigma(G_1i,G_1j)$ and $\sigma(G_2i,G_2j)$ are connected by
a zig-zag of direct equivalences.
\end{cor}

\begin{proof}
Call the two enrichments $\sigma_1$ and $\sigma_2$. Define a
$\sigma_2-\sigma_1$ bimodule by $M(i,j)=\sigma(G_1i,G_2j)$.  The maps
$G_1i \we G_2i$ give rise to maps $S\ra M(i,i)$, making $M$ into a
pointed bimodule.  One readily checks that this is a quasi-equivalence
between $\sigma_2$ and $\sigma_1$, and then applies
Proposition~\ref{pr:equiv=direct}.
\end{proof}

\subsection{Homotopy invariant enrichments}
We give a few other basic results about model enrichments.

\begin{prop}
\label{pr:coffib}
Let $Qa \trfib a$ be a cofibrant-replacement functor in $\cM$, and let
$x \trcof Fx$ be a fibrant-replacement functor.  If $\tau$ is a model
enrichment of $\cM$, then $\tau(Qa,Qb)$ and $\tau(Fa,Fb)$ give model
enrichments which are quasi-equivalent to $\tau$.
\end{prop}

\begin{proof}
Left to the reader (see Example~\ref{ex:example}).
\end{proof}

A model enrichment $\tau$ of $\cM$ by $\cV$ will be called
\dfn{homotopy invariant} if whenever $a\ra a'$ and $x\ra x'$ are weak
equivalences then the maps $\tau(a',x)\ra \tau(a,x)\ra \tau(a,x')$ are
both weak equivalences as well.  Note that there is no
cofibrancy/fibrancy assumption on the objects.

\begin{cor}
\label{co:hom-enrich}
Every model enrichment is quasi-equivalent to one which is homotopy
invariant.
\end{cor}

\begin{proof}
Let $\tau$ be a model enrichment of $\cM$ by $\cV$.  By
Proposition~\ref{pr:coffib} (used twice), the enrichments $\tau(a,b)$,
$\tau(Qa,Qb)$, and $\tau(QFa,QFb)$ are all quasi-equivalent.  The last
of these is homotopy invariant.
\end{proof}

Recall that the monoidal product on $\ho(\cV)$ is defined by $v_1
\tens_{\mathbb L} v_2=Cv_1 \tens Cv_2$, where $C$ is some chosen
cofibrant-replacement functor in $\cV$.  It is easy to check that a
homotopy invariant enrichment $\tau$ induces an enrichment of
$\ho(\cM)$ by $\ho(\cV)$, where the composition maps are the
composites
\[ \tau(b,c) \tens_{\mathbb L} \tau(a,b) \ra \tau(b,c)\tens \tau(a,b)
\ra \tau(a,c).
\]
We note the following:

\begin{cor}
\label{co:hoenrich2}
If two homotopy invariant enrichments $\sigma$ and $\tau$ are
quasi-equivalent, then the induced enrichments of $\ho(\cM)$ by
$\ho(\cV)$ are equivalent.  
\end{cor}

\begin{proof}
First note that if $M$ is a quasi-equivalence between $\sigma$ and
$\tau$ then $M$ is automatically homotopy invariant itself (in the
obvious sense)---this follows from the two-out-of-three property for
weak equivalences. Therefore $M$ may be extended to a functor on the
homotopy category, where it clearly gives an equivalence between
the enrichments induced by $\sigma$ and $\tau$.

To say that $\sigma$ and $\tau$ are quasi-equivalent, though, does not
say that such an $M$ necessary exists---it only says that there is a
chain of such $M$'s.  Note that the intermediate model enrichments in
the chain need not be homotopy invariant.  To get around this, we do
the following.  If $\mu$ is a model enrichment of $\cM$ by $\cV$, let
$\mu^h$ be the model enrichment $\mu^h(a,b)=\mu(QFa,QFb)$.  We have
seen that this is homotopy invariant and quasi-equivalent to $\mu$.
If $M$ is a quasi-equivalence between $\mu_1$ and $\mu_2$, note that
$M^h$ (with the obvious definition) is a quasi-equivalence between
$\mu_1^h$ and $\mu_2^h$.  It follows readily that if our $\sigma$ and
$\tau$ are quasi-equivalent then they are actually quasi-equivalent
through a chain where all the intermediate steps are homotopy
invariant.  Now one applies the first paragraph to all the links in
this chain.
\end{proof}

\subsection{Transporting enrichments}

Let $G\colon \cM \ra \cN$ be a functor, and suppose $\tau$ is an
enrichment of $\cN$.  Define an enrichment $G^*\tau$ of $\cM$ by the
formula $G^*\tau(m_1,m_2)=\tau(Gm_1,Gm_2)$.  Call this the
\dfn{pullback} of $\tau$ along $G$.

\begin{lemma}
\label{le:transp1}
Let $\cM$ and $\cN$ be model categories, and let $G\colon \cM \ra \cN$
be a functor which preserves weak equivalences and has its image in
the cofibrant-fibrant objects of $\cN$.  If $\tau$ is a model
enrichment of $\cN$, then $G^*\tau$ is a model enrichment of $\cM$.
Moreover, $G^*$ preserves quasi-equivalence: it induces $G^*\colon
\ME_0(\cN,\cV) \ra \ME_0(\cM,\cV)$.  
\end{lemma}

\begin{proof}
Routine.
\end{proof}

\begin{lemma}
\label{le:transp2}
Let $\cM$ and $\cN$ be model categories, and let $\tau$ be a homotopy
invariant enrichment of $\cN$.  Suppose $G_1,G_2 \colon \cM \ra \cN$
are two functors which preserve weak equivalences, and assume there is
a natural weak equivalence $G_1 \we G_2$.  Then $G_1^*\tau$ and
$G_2^*\tau$ are model enrichments of $\cM$, and they are
quasi-equivalent.
\end{lemma}

\begin{proof}
The quasi-equivalence is given by $M(a,b)=\tau(G_1a,G_2b)$.  The weak
equivalences $G_1a \ra G_2a$ give the necessary maps $S\ra M(a,a)$.  
Details are left to the reader.  
\end{proof}

Recall that a \dfn{Quillen map} $L\colon \cM \ra \cN$ is an adjoint
pair $L\colon \cM\adjoint \cN \colon R$ in which $L$ preserves
cofibrations and trivial cofibrations (and  $R$ preserves
fibrations and trivial fibrations).  
Choose cofibrant-replacement functors $Q_\cM X
\trfib X$ and $Q_\cN Z\trfib Z$ as well as fibrant-replacement
functors $A\trcof F_\cM A$ and $B\trcof F_\cN B$.  If $\tau$ is a
model enrichment of $\cN$ by $\cV$, we can define a model enrichment
on $\cM$ by the formula $L^*\tau(a,x)=\tau(F_\cN L Q_\cM a,F_\cN
LQ_\cM x)$.  Similarly, if $\sigma$ is a model enrichment of $\cM$ by
$\cV$ we get a model enrichment on $\cN$ by the formula
$L_*\sigma(c,w)=\sigma(Q_\cM R F_\cN c, Q_\cM RF_\cN w)$.

\begin{prop}
\label{pr:transp}
\mbox{}\par
\begin{enumerate}[(a)]
\item
The constructions $L^*$ and $L_*$ induce maps $L^*\colon
\ME_0(\cN,\cV) \ra \ME_0(\cM,\cV)$ and $L_*\colon \ME_0(\cM,\cV) \ra
\ME_0(\cN,\cV)$.
\item 
The maps in (a) do not depend on the choice of cofibrant- and fibrant-
replacement functors.
\item
If $L,L'\colon \cM \ra \cN$ are two maps which are
Quillen-homotopic, then $L_*=L'_*$ and $L^*=(L')^*$ 
as maps on $\ME_0(\blank,\cV)$.
\item If $L\colon \cM \ra \cN$ is a Quillen equivalence, then the
functors $L^*$ and $L_*$ are inverse isomorphisms $\ME_0(\cM,\cV) \iso
\ME_0(\cN,\cV)$.
\item Suppose $\cM$ and $\cN$ are $\cV$-model categories, with the
associated $\cV$-enrichments denoted $\sigma_\cM$ and $\sigma_\cN$.  
If $L\colon \cM \ra \cN$ is a $\cV$-Quillen equivalence, then
$L_*(\sigma_\cM)=\sigma_\cN$ and $L^*(\sigma_\cN)=\sigma_\cM$ (as
elements of $\ME_0(\blank,\cV)$).
\end{enumerate}
\end{prop}

For the notion of `$\cV$-Quillen equivalence' used in part (e), see
Section~\ref{se:Afunctor}.

\begin{proof}
We will only prove the results for $L^*$; proofs for $L_*$ are
entirely similar.

Part (a) follows from Lemma~\ref{le:transp1}, as the composite functor
$F_\cN L Q_\cM$ preserves weak equivalences and has its image in the
cofibrant-fibrant objects.

For part (b), suppose $Q_1X\trfib X$ and $Q_2 X\trfib X$ are two
cofibrant-replacement functors for $\cM$.  Write $L_1^*$ and $L_2^*$
for the resulting maps $\ME_0(\cN,\cV) \ra \ME_0(\cM,\cV)$.  By
Corollary~\ref{co:hom-enrich} it suffices to show that
$L_1^*(\tau)=L_2^*(\tau)$ for any homotopy invariant enrichment
$\tau$.  Let $Q_3X=Q_1X \times_X Q_2X$.  There is a zig-zag of natural
weak equivalences $Q_1 \bwe Q_3 \we Q_2$.  The result now follows by
Lemma~\ref{le:transp2} applied to the composites $FLQ_1$, $FLQ_3$, and
$FLQ_2$.

For part (c), it again suffices to prove $L^*(\tau)=(L')^*(\tau)$ in
the case where $\tau$ is homotopy invariant.  The Quillen homotopy is
a natural transformation $L \ra L'$ which is a weak equivalence on
cofibrant objects.  The result is then a direct application of
Lemma~\ref{le:transp2}.

For (d) we will check that if $\tau$ is a homotopy invariant enrichment
of $\cN$ then $L_*(L^*\tau)=\tau$ in $\ME_0(\cN,\cV)$.  The enrichment
$L_*(L^*\tau)$ is the pullback of $\tau$ along the composite functor
$FLQQRF\colon \cN \ra \cN$. There is a zig-zag of natural weak
equivalences
\[ FLQQRF \bwe LQQRF \we F \bwe \Id 
\]
(the second being the composite $LQQRF \ra LRF \ra F$, which is a weak
equivalence because we have a Quillen equivalence).  Each of the
functors in the zig-zag preserves weak equivalences, so the result
follows from Lemma~\ref{le:transp2}.  

Finally, we prove (e).  By (d) it suffices just to prove
$L_*\sigma_\cM=\sigma_\cN$.  The assumption gives us a natural
isomorphism $\sigma_\cN(LA,X)\iso \sigma_\cM(A,RX)$ (see
Section~\ref{se:Afunctor}).  One checks that the enrichments
$\sigma_\cM(QRFX,QRFY)$ and $\sigma_\cN(FX,FY)$ are quasi-equivalent
via the bimodule $M(X,Y)=\sigma_\cM(QRFX,RFY)\iso
\sigma_\cN(LQRFX,FY)$. (The verification that this really {\it is\/} a
bimodule requires some routine but tedious work, mainly using
Remark~\ref{re:adjunct}).  But Proposition~\ref{pr:coffib} says that
$\sigma_\cN(FX,FY)$ is quasi-equivalent to $\sigma_\cN$, so we are
done.
\end{proof}


\section{Universal pointed model categories}
\label{se:univ}

If $\cC$ is a small category then there is a `universal model
category' built from $\cC$.  This was developed in \cite{D1}.  The
present section deals with a pointed version of that theory.  The
category of functors from $\cC$ to pointed simplicial sets plays the
role of a universal {\it pointed\/} model category built from $\cC$.

\medskip

\subsection{Basic definitions}
Recall from \cite{D1} that if $\cC$ is a small category then $\U\cC$
denotes the model category of simplicial presheaves on $\cC$, with
fibrations and weak equivalences defined objectwise.  One has the
Yoneda embedding $r\colon \cC \inc \U\cC$ where $rX$ is the presheaf
$Y\mapsto \cC(Y,X)$.  

Let $\U_+\cC$ be the category of functors from $\cC^{op}$ into {\it
pointed\/} simplicial sets, with the model structure where weak
equivalences and fibrations are again objectwise.  This can also be
regarded as the undercategory $(*\ovcat \U\cC)$.   

There is a Quillen map $\U\cC \ra \U_+\cC$ where the left adjoint
sends $F$ to $F_+$ (adding a disjoint basepoint) and the right adjoint
forgets the basepoint.  Write $r_+$ for the composite $\cC \inc \U\cC
\ra \U_+\cC$.

Finally, if $S$ is a set of maps in $\U\cC$ then let $S_+$ denote the
image of $S$ under $\U\cC \ra \U_+\cC$.  Note that if all the maps in
$S$ have cofibrant domain and codomain, then by
\cite[Prop. 3.3.18]{H} one has an induced Quillen map $\U\cC/S \ra
\U_+\cC/(S_+)$.

The following simple lemma unfortunately has a long proof:

\begin{lemma}
\label{le:ugly}
Let $S$ be a set of maps between cofibrant objects in $\U\cC$, and
suppose that the map $\emptyset \ra *$ is a weak equivalence in
$\U\cC/S$.  Then $\U\cC/S \ra \U_+\cC/(S_+)$ is a Quillen equivalence.
\end{lemma}

\begin{proof}
Write $\cM=\U\cC$ and $\cM_+=\U_+\cC=(*\ovcat \cM)$ (the lemma
actually holds for any simplicial, left proper, cellular model
category in place of $\U\cC$).  Write $F\colon \cM \adjoint
\cM_+\colon U$ for the Quillen functors.  We will start by showing
that a map in $\cM_+/(S_+)$ is a weak equivalence if and only if it's
a weak equivalence in $\cM/S$.  Unfortunately the proof of this fact
is somewhat lengthy.

An object $X \in \cM_+$ is $(S_+)$-fibrant if it is fibrant in $\cM_+$
(equivalently, fibrant in $\cM$) and if the induced map on simplicial
mapping spaces $\mM_{\cM_+}(B_+,X) \ra \mM_{\cM_+}(A_+,X)$ is a weak
equivalence for every $A\ra B$ in $S$.  By adjointness, however,
$\mM_{\cM_+}(A_+,X)\iso \mM_{\cM}(A,X)$ (and similarly for $B$).
It follows that $X\in \cM_+$ is $(S_+)$-fibrant if and only if $X$ is
$S$-fibrant in $\cM$.

Suppose $C$ is a cofibrant object in $\cM$.  Using the fact that
$\cM/S$ is left proper and that $\emptyset \ra *$ is a weak
equivalence, it follows that $C \ra C\amalg *$ is also a weak
equivalence in $\cM/S$.  As a consequence, if $C\ra D$ is a map
between cofibrant objects which is a weak equivalence in $\cM/S$, then
$C_+ \ra D_+$ is also a weak equivalence in $\cM/S$.

Now consider the construction of the localization functor $L_{S_+}$
for $\cM_+/(S_+)$.  This is obtained via the small object argument, by
iteratively forming pushouts along the maps 
\[ [\Lambda^{n,k} \ra \Delta^n] \tens_+ [A_+ \ra B_+].
\]
Here ``$\tens_+$'' denotes the simplicial tensor in the pointed
category $\cM_+$, that is to say $K\tens_+ A = (K_+ \tens A)/((*\tens A)
\amalg (K_+\tens *))$ for $K\in \sSet$ and $A\in \cM$.  
The above maps are then readily identified with the maps
\[  \Bigl [(\Lambda^{n,k} \tens B) \amalg_{\Lambda^{n,k}\tens A} (\Delta^n
\tens A) \Bigr ]_+ \ra (\Delta^n \tens B)_+.
\]
As 
$ [(\Lambda^{n,k} \tens B) \amalg_{\Lambda^{n,k}\tens A} (\Delta^n
\tens A) \Bigr ] \ra (\Delta^n \tens B)$ is a map between cofibrant
objects which is a weak equivalence in $\cM/S$, so is the displayed
map above.  It follows that for any $X\in \cM_+$, the map $X \ra
L_{S_+}X$ is a weak equivalence in $\cM/S$ (in addition to being a
weak equivalence in $\cM_+/(S_+)$, by construction).  

Let $X\ra Y$ be a map in $\cM_+$.  Consider the square
\[ \xymatrix{ X \ar[r] \ar[d] & Y \ar[d] \\
     L_{S_+}X \ar[r] & L_{S_+} Y.
}
\]
The vertical maps are weak equivalences in both $\cM/S$ and
$\cM_+/(S_+)$.  If $X\ra Y$ is a weak equivalence in $\cM_+/(S_+)$, then
the bottom map is a weak equivalence in $\cM_+$.  This is the same as
being a weak equivalence in $\cM$, and therefore $X\ra Y$ is also a
weak equivalence in $\cM/S$ (going back around the square, using the
2-out-of-3 property).  Similarly, if $X\ra Y$ is a weak equivalence in
$\cM/S$ then so is the bottom map.  But the objects $L_{S_+}X$ and
$L_{S_+}Y$ are fibrant in $\cM/S$, so the bottom map is actually a
weak equivalence in $\cM$ (and also in $\cM_+$).  It follows that
$X\ra Y$ is a weak equivalence in $\cM_+/(S_+)$.  

This completes the proof that a map in $\cM_+/(S_+)$ is a weak
equivalence if and only if it is so in $\cM/S$.

To show that $\cM/S \ra \cM_+/(S_+)$ is a Quillen equivalence we must
show two things.  If $A$ is
a cofibrant object in $\cM$ and $A_+ \ra X$ is a fibrant replacement
in $\cM_+/(S_+)$, we must show that $A\ra X$ is a weak equivalence in
$\cM/S$.  But from what we have already shown we know $A\ra A_+$ and
$A_+\ra X$ are weak equivalences in $\cM/S$, so this is obvious.
We must also show that if $Z$ is a fibrant object in $\cM_+/(S_+)$
and $B\ra Z$ is a cofibrant replacement in $\cM/S$, then $B_+ \ra Z$
is a weak equivalence in $\cM_+/(S_+)$.   This is the same as showing
it's a weak equivalence in $\cM/S$.  But in the sequence $B \ra
B_+ \ra Z$, the first map and the composite are both equivalences in
$\cM/S$; so the map $B_+ \ra Z$ is an equivalence as well. 
\end{proof}

\subsection{Basic properties}

\begin{prop}
\label{pr:factor}
Suppose that $L\colon \U\cC/S \ra \cM$ is a Quillen map, where $S$ is
a set of maps between cofibrant objects.  If $\cM$ is pointed, there
is a Quillen map $L_+\colon \U_+\cC/(S_+) \ra \cM$ such that the
composite $\U\cC/S \ra \U_+\cC/(S_+) \ra \cM$ is $L$.  If $L$ is a
Quillen equivalence, then so is $L_+$.
\end{prop}

\begin{proof}
For any $A\in \cC$, write $\Gamma_*A$ for the cosimplicial object
$[n]\mapsto L(rA \tens \Delta^n)$.  Recall that the right adjoint to
$L$ sends an $X\in \cM$ to the simplicial presheaf $A\mapsto
\cM(\Gamma_*A,X)$.  Since $\cM$ is pointed, this simplicial presheaf
is also pointed.  Let $\sing_* \colon \cM \ra \U_+\cC$ be this functor.

If $F\in \U_+\cC$ define $L_+(F)$ to be the pushout of $* \la L(*) \ra
L(F)$.  This is readily seen to be left adjoint to $\sing_*$.  It is
also easy to check that $L_+\colon \U_+\cC \ra \cM$ is a Quillen map
and the composite $\U\cC \ra \U_+\cC \ra \cM$ equals $L$.

To obtain the map $\U_+\cC/(S_+) \ra \cM$ one only has to see that
$L_+$ maps elements of $S_+$ to weak equivalences in $\cM$.  But this
is obvious: if $A\in \U\cC$ then $L_+(A\amalg *)\iso L(A)$, and $L$
takes elements of $S$ to weak equivalences.

Finally, assume that $L$ is a Quillen equivalence.  Since $\cM$ is
pointed, it follows that $\emptyset \ra *$ is a weak equivalence in
$\U\cC/S$ (using that $L(\emptyset)=*$ and $R(*)=*$).  So by the above
lemma, $\U\cC/S \ra \U_+\cC/(S_+)$ is a Quillen equivalence; therefore
$L_+$ is one as well.
\end{proof}

The next two propositions of this section accentuate the roll of
$\U_+\cC$ as the universal pointed model category built from $\cC$.
These results are direct generalizations of \cite[Props. 2.3,
5.10]{D1}.

\begin{prop}
\label{pr:univ2}
Let $\cC$ be a small category, and let $\gamma\colon \cC \ra\cM$ be a functor
from $\cC$ into a pointed model category $\cM$.  Then $\gamma$
``factors'' through $\cU_+\cC$, in the sense that there is a
Quillen pair $L \colon \U_+\cC \adjoint \cM \colon R$ and a natural
weak equivalence $L\circ r_+ \we \gamma$.  Moreover, the category of
all such factorizations---as defined in \cite[p. 147]{D1}---is contractible.
\end{prop}

\begin{proof}
The result follows from \cite[Prop. 2.3]{D1} and
Proposition~\ref{pr:factor} above.
\end{proof}

\begin{prop}
\label{pr:univ3}
Suppose $L\colon \U_+\cC/S \ra \cN$ is a Quillen map, and $P\colon \cM
\we \cN$ is a Quillen equivalence between pointed model categories.
Then there is a Quillen map $L'\colon \U_+\cC/S \ra \cM$ such that
$P\circ L'$ is Quillen homotopic to $L$.  Moreover, if $\cM$ is
simplicial then $L'$ can be chosen to be simplicial.
\end{prop}

\begin{proof}
The first statement follows directly from \cite[Prop. 5.10]{D1} and
Proposition~\ref{pr:factor} above.  The second statement was never
made explicit in \cite{D1}, but follows at once from analyzing the
proof of \cite[Prop. 2.3]{D1}.  To define $F$ one first gets a map
$f\colon \cC \ra \cM$ with values in the cofibrant objects, and then
$F$ can be taken to be the unique colimit-preserving functor
characterized by $F(rA\tens K)=f(A)\tens K$, where $A\in \cC$ and
$K\in \sSet$.  This is clearly a simplicial functor.
\end{proof}

\begin{prop} 
\label{pr:univ}
Let $\cM$ be a pointed, combinatorial model category.
\begin{enumerate}[(a)]
\item
There is a
Quillen equivalence $\U_+\cC/S \ra \cM$ for some $\cC$ and $S$.
\item Let $\cN$ be a pointed model category, and let $\cM \bwe
\cM_1 \we  \cdots \bwe \cM_n \we \cN$ be a zig-zag of
Quillen equivalences (where the intermediate model categories are not
necessarily pointed or combinatorial).  Then there is a simple zig-zag
of Quillen equivalences
\[ \cM \bwe \U_+ \cC/S \we \cN 
\]
for some $\cC$ and $S$.
\item In the context of (b), the simple zig-zag can be chosen so that
the derived equivalence $\ho(\cM) \he \ho(\cN)$ is isomorphic to the
derived equivalence specified by the original zig-zag.  
\end{enumerate}
\end{prop}

In part (b), note that we have replaced a zig-zag of Quillen
equivalences---in which the intermediate steps are not necessarily
pointed---by one in which the intermediate steps {\it are\/} pointed.
For (c), recall that two pairs of adjoint functors $L\colon \cC
\adjoint \cD\colon R$ and $L'\colon \cC \adjoint \cD \colon R'$ are
said to be isomorphic if there is a natural isomorphism $LX \iso L'X$
for all $X\in \cC$ (equivalently, if there is a natural isomorphism
$RY \iso R'Y$ for all $Y\in \cD$).

\begin{proof}
Let $\cM$ be a pointed, combinatorial model category.  By
\cite[Th. 6.3]{D1} there is a Quillen equivalence $\U\cC/S \ra \cM$
for some $\cC$ and $S$.  Proposition~\ref{pr:factor} shows there is an
induced Quillen equivalence $\U_+\cC/(S_+) \ra \cM$.  This proves (a).

Parts (b) and (c) follow in the same way from \cite[Cor. 6.5]{D1}, or
directly by applying Proposition~\ref{pr:univ3}.
\end{proof}

\subsection{Application to stabilization}

Suppose $\cM$ is a stable model category, and we happen to have a
Quillen equivalence $\U_+\cC/S \ra \cM$.  It follows in particular
that $\U_+\cC/S$ is also stable.  Now, $\U_+\cC/S$ is a simplicial,
left proper, cellular model category.  So using \cite[Secs. 7,8]{Ho1}
we can form the corresponding category of symmetric spectra
$\Spe^\Sigma(\U_+\cC/S)$ (with its stable model structure).  This comes
with a Quillen map $\U_+\cC/S \ra \Spe^\Sigma(\U_+\cC/S)$, and since
$\U_+\cC$ is stable this map is a Quillen equivalence
\cite[Th. 9.1]{Ho1}.  Finally, the category $\U_+\cC/S$ satisfies the
hypotheses of \cite[Th. 8.11]{Ho1}, and so $\Spe^\Sigma(\U_+\cC/S)$ is
a {\it spectral\/} model category (in the sense of
Section~\ref{se:spectral}).  We have just proven part (a) of the
following:
  
\begin{prop}
\label{pr:presentation}
Let $\cM$ be a stable model category, and suppose $\U_+\cC/S \ra \cM$
is a Quillen equivalence.
\begin{enumerate}[(a)]
\item There is a zig-zag of Quillen equivalences $\cM \bwe \U_+\cC/S
\we \Spe^\Sigma(\U_+\cC/S)$.
\item If $\U_+\cD/T \ra \cM$ is another Quillen equivalence, there is
a diagram of Quillen equivalences
\[\xymatrix{ 
\cM  & \U_+\cC/S \ar[l]\ar[r]\ar[d] & \Spe^\Sigma(\U_+\cC/S) \ar[d] \\
 & \U_+\cD/T \ar[ul]\ar[r] & \Spe^\Sigma(\U_+\cD/T)
}
\]
where the left vertical map is a simplicial adjunction, the right
vertical map is a spectral adjunction, the square commutes
on-the-nose, and the triangle commutes up to a Quillen homotopy.
\end{enumerate}
\end{prop}

\begin{proof}
We have left only to prove (b).  Given Quillen equivalences $L_1\colon
\U_+\cC/S \ra \cM$ and $L_2\colon \U_+\cD/T \ra \cM$, it follows from
Proposition~\ref{pr:univ3} that there is a Quillen map $F\colon
\U_+\cC/S \ra \U_+\cD/T$ making the triangle commute up to Quillen
homotopy.  Since $\U_+\cD/T$ is a simplicial model category, we can
choose $F$ to be simplicial.
But this ensures that $\Spe^\Sigma(\U_+\cC/S) \ra \Spe^\Sigma(\U_+\cD/T)$ is
spectral.
\end{proof}


\section{The main results}
\label{se:main}
In this section we attach to any stable, combinatorial model category
$\cM$ a model enrichment $\tau_\cM$ over symmetric spectra.  This involves
choices, but these choices only affect the end result up to
quasi-equivalence.  We also show that a zig-zag of Quillen
equivalences between model categories $\cM$ and $\cN$ must carry
$\tau_\cM$ to $\tau_\cN$.  So the canonical enrichments $\tau$
give rise to invariants of model categories up to Quillen equivalence.
Finally, we specialize all these results to establish basic properties
of homotopy endomorphism spectra.

The present results are all direct consequences of work from previous
sections.  Our only job is to tie everything together.

\medskip

\subsection{Construction of spectral enrichments}
Let $\cM$ be a stable, combinatorial model category.  By
Proposition~\ref{pr:presentation}(a) there is a zig-zag of Quillen
equivalences
\[ \cM \llla{L} \U_+\cC/S \llra{F} \Spe^\Sigma(\U_+\cC/S).
\]
The right-most model category comes equipped with a spectral enrichment
$\sigma$.   We define $\tau_\cM \in \ME_0(\cM,\Spe^\Sigma)$ to be
$L_*(F^*\sigma)$.  

\begin{prop}
The element $\tau_\cM \in \ME_0(\cM,\Spe^\Sigma)$ doesn't depend on the
choice of $\cC$, $S$, or the Quillen equivalence $\U_+\cC/S \we \cM$.
\end{prop}

\begin{proof}
Applying $\ME_0(\blank,\Spe^\Sigma)$ to the diagram from
Proposition~\ref{pr:presentation}(b) gives a commutative diagram of
bijections, by Proposition~\ref{pr:transp}.  The result follows
immediately from chasing around this diagram and using
Proposition~\ref{pr:transp}(e). 
\end{proof}

Choose a homotopy invariant enrichment quasi-equivalent to $\tau_\cM$.
By Corollary~\ref{co:hoenrich2} this induces an enrichment of
$\ho(\cM)$ by $\ho(\Spe^\Sigma)$, and different choices lead to
equivalent enrichments.  This proves Corollary~\ref{co:hoenrich}.

We now turn our attention to functoriality:

\begin{prop}
\label{pr:quillenequiv}
Suppose $L\colon \cM \ra \cN$ is a Quillen equivalence between stable,
combinatorial model categories.  Then $L^*(\tau_\cN)=\tau_\cM$ and
$L_*(\tau_\cM)=\tau_\cN$.
\end{prop}

\begin{proof}
Choose a Quillen equivalence $\U_+\cC/S \ra \cM$, by
Proposition~\ref{pr:univ2}.  
We then have a diagram of Quillen equivalences
\[ \xymatrix{
\Spe^\Sigma(\U_+\cC/S) & \U_+\cC/S\ar[r]\ar[l] & \cM\ar[r] & \cN.
}
\]
Applying $\ME_0(\blank,\Spe^\Sigma)$ to the
diagram yields a diagram of bijections by
Proposition~\ref{pr:transp}.  The result follows from chasing
around this diagram.
\end{proof}

\begin{remark}
The above result is more useful in light of
Proposition~\ref{pr:univ}(b).  Suppose $\cM$ and $\cN$ are stable,
combinatorial model categories which are Quillen equivalent.  This
includes the possibility that the Quillen equivalence occurs through a
zig-zag, where the intermediate steps may not be combinatorial or
pointed.  So the above result doesn't apply directly.  However,
Proposition~\ref{pr:univ}(b) shows that any such zig-zag may be
replaced by a simple zig-zag where the intermediate step is both
combinatorial and pointed (hence also stable).  One example of this
technique is given in the proof of Theorem~\ref{th:invariance} below.
\end{remark}

\begin{prop}
\label{pr:spectralenrich}
Assume that $\cM$ is stable, combinatorial, and a spectral model
category.  Then $\tau_\cM$ is quasi-equivalent to the enrichment
$\sigma$ provided by the spectral structure.
\end{prop}

\begin{proof}
As $\cM$ is spectral, it is in particular simplicial
(cf.~\ref{se:spectral}).  So one may choose a Quillen equivalence
$L\colon \U_+\cC/S \ra \cM$ consisting of simplicial functors (see
discussion in the proof of Proposition~\ref{pr:presentation}).  We have the
Quillen maps
\[ \xymatrix{ \U_+\cC/S \ar[r]\ar[d] & \cM \\
       \Spe^\Sigma(\U_+\cC/S).
}
\]
We claim there is a spectral Quillen equivalence
$\Spe^\Sigma(\U_+\cC/S) \ra \cM$ making the triangle commute. 
This immediately implies the result we want: applying
$\ME_0(\blank,\Spe^\Sigma)$ to the triangle gives a commutative diagram
of bijections by Proposition~\ref{pr:transp}(d), and the diagonal map
sends the canonical spectral enrichment of $\Spe^\Sigma(\U_+\cC/S)$ to
the given spectral enrichment of $\cM$ by Proposition~\ref{pr:transp}(e).

We are reduced to constructing the spectral Quillen map
$\Spe^\Sigma(\U_+\cC/S) \ra \cM$.  Note that objects in
$\Spe^\Sigma(\U_+\cC)$ may be regarded as presheaves of symmetric
spectra on $\cC$.  That is, we are looking at the functor category
$\Func(\cC^{op}, \Spe^\Sigma)$.  By
Proposition~\ref{pr:diagram-adjunct}, the composite $\cC \ra \U_+\cC
\ra \cM$ induces a spectral Quillen map $\Real \colon
\Func(\cC^{op},\Spe^\Sigma) \adjoint \cM \colon \Sing$, where the
functor category is given the `objectwise' model structure.  Note that
the composite of right adjoints $\cM \ra \Spe^\Sigma(\U_+\cC/S) \ra
\U_+\cC/S$ is indeed the right adjoint of $L$.

We need to check that $(\Real,\Sing)$ give a Quillen map
$\Spe^\Sigma(\U_+\cC/S) \ra \cM$.  By Proposition~\ref{pr:locstab}, the
domain model category is identical to $(\Spe^\Sigma \U_+\cC)/S_{stab}$
(notation as in Appendix B).  But to show a Quillen map $\Spe^\Sigma(
\U_+\cC) \ra \cM$ descends to $(\Spe^\Sigma \U_+\cC)/S_{stab}$, it is
sufficient to check that the left adjoint sends elements of $S_{stab}$
to weak equivalences in $\cM$.

A typical element of $S_{stab}$ is a map $F_i(A) \ra F_i(B)$ where
$A\ra B$ is in $S$ ($F_i(\blank)$ is defined in Appendix B).  Certainly
$\Real$ sends $F_0A \ra F_0B$ to a weak equivalence, since $\Real
\circ F_0$ is the map $L\colon \U_+\cC \ra \cM$ and this map sends
elements of $S$ to weak equivalences by construction.  For $i\geq 1$,
note that the $i$th suspension of $F_iA \ra F_iB$ is $F_0 A\ra F_0 B$.
Since $\cM$ is a stable model category, the fact that $\Real$ sends
$F_0 A \ra F_0 B$ to a weak equivalence therefore immediately implies
that it does the same for $F_i A \ra F_i B$.
\end{proof}

\subsection{Homotopy endomorphism spectra}
\label{se:endo}
Let $\cM$ be a stable, combinatorial model category, and let $X\in
\cM$ be a cofibrant-fibrant object.  Consider the ring spectrum
$\tau_\cM(X,X)$.  By Corollary~\ref{co:endo=invariant}, the
isomorphism class of this ring spectrum in $\Ho(\RingSpectra)$ only
depends on $\tau_\cM$ up to quasi-equivalence.

Now let $W$ be an arbitrary object in $\cM$, and let $X_1$ and $X_2$
be two cofibrant-fibrant objects weakly equivalent to $W$.  Then there
exists a weak equivalence $f\colon X_1 \ra X_2$.  Let $\cI$ be the
category with one object and an identity map, and consider the two
functors $\cI \ra \cM$ whose images are $X_1$ and $X_2$, respectively.
Applying Corollary~\ref{co:endo=coffib} to this situation, we find
that $\tau_\cM(X_1,X_1)$ and $\tau_\cM(X_2,X_2)$ are weakly equivalent
ring spectra.  So the corresponding isomorphism class in
$\Ho(\RingSpectra)$ is a well-defined invariant of $W$.  We will write
$\hEnd(W)$ for any ring spectrum in this isomorphism class.

The two main results about homotopy endomorphism ring spectra were
stated as Theorem~\ref{th:invariance} and
Proposition~\ref{pr:spectralcase}.  We now give the proofs:

\begin{proof}[Proof of Theorem~\ref{th:invariance}]
If two stable, combinatorial model categories $\cM$ and $\cN$ are
Quillen equivalent through a zig-zag, then by
Proposition~\ref{pr:univ}(b,c) there is a simple zig-zag
$\cM \bwe \U_+\cC/S \we \cN$
inducing an isomorphic derived equivalence of the homotopy categories.   
Now we apply Proposition~\ref{pr:quillenequiv} (twice) to connect
$\tau_\cM$ to $\tau_\cN$.  Finally, the required equivalence of
homotopy endomorphism ring spectra follows from
Corollary~\ref{co:endo=invariant}.    
\end{proof}

\begin{proof}[Proof of Proposition~\ref{pr:spectralcase}]
This is a special case of Proposition~\ref{pr:spectralenrich}.
\end{proof}

\vspace{0.1in}


\section{A leftover proof}
\label{se:leftover}

In this section we complete the proof of
Proposition~\ref{pr:equiv=direct}.  Essentially this amounts to just
explaining why the proof has already been given in \cite[Lemma
A.2.3]{SS2}.  The differences between our situation and that of
\cite{SS2} are (1) our indexing categories are not necessarily
discrete (i.e., they have maps other than identities), and (2) we are
dealing with a general symmetric monoidal model category rather than
symmetric spectra.  It turns out that neither difference is
significant.

\medskip

\subsection{Modules}
Let $\cV$ be a symmetric monoidal category.  Let $\cC$ be a category,
and let $\sigma$ be an enrichment of $\cC$ by $\cV$.  A \mdfn{left
$\sigma$-module} is a collection of objects $M(c)\in \cV$ (for each
$c\in \cC$) together with maps $\sigma(a,b) \tens M(a) \ra M(b)$ such
that the following diagrams commute:
\[\xymatrixcolsep{1.5pc}
\xymatrix{ \sigma(b,c) \tens \sigma(a,b) \tens M(a) \ar[r] \ar[d] &
                \sigma(b,c)\tens M(b) \ar[d]  & S\tens M(a)
                \ar[r]\ar[dr] & \sigma(a,a)\tens M(a) \ar[d] \\
 \sigma(a,c)\tens M(a) \ar[r] & M(c) && M(a)
}
\]
As for the case of bimodules (see Section~\ref{se:bimodules}), $M$
inherits a natural structure of a functor $\cC \ra \cV$.
(An $\S_\cC$-module is precisely a functor $M\colon\cC \ra \cV$, and
so the map $\S_\cC \ra \sigma$ gives every left $\sigma$-module a
structure of functor by restriction).

\begin{remark}
A more concise way to phrase the above definition is to say that a
left $\sigma$-module is a $\cV$-functor from the $\cV$-enriched category
$\cC$ to the $\cV$-enriched category $\cV$.
\end{remark}

We now record several basic facts about modules and functors.  To
begin with, one can check that colimits and limits in the category of
$\sigma$-modules are the same as those in the category of functors
$\Func(\cC,\cV)$.  

For each $c\in \cC$, note that the functor $\sigma(c,\blank)\colon \cC
\ra \cV$ has an obvious structure of left $\sigma$-module.  It is the
`free' module determined by $c$.  For $A\in \cV$ we write
$\sigma(c,\blank) \tens A$ for the module $a\mapsto \sigma(c,a)\tens
A$.

The canonical map $\S_\cC \ra\sigma$ induces a forgetful functor from
$\sigma$-modules to $\S_\cC$-modules, which is readily checked to have
a left adjoint: we'll call this adjoint $\sigma\tens(\blank)$.  Let
$T\colon (\S_\cC-mod) \ra (\S_\cC-mod)$ be the resulting cotriple.
It's useful to note that if $M\colon \cC \ra \cV$ is a functor then
$\sigma\tens M$ is the coequalizer of
\[ \coprod_{a\ra b} \sigma(b,\blank)\tens M(a) \dbra \coprod_a
\sigma(a,\blank)\tens M(a)
\]
(the coequalizer can be interpreted either in the category of
$\sigma$-modules or the category of functors, as they coincide).

Given two functors $M,N\colon \cC \ra \cV$, one can define $\F(M,N)
\in \cV$ as the equalizer of $\prod_a \cV(M(a),N(a)) \dbra \prod_{a\ra
b} \cV(M(a),N(b))$.  Together with the `objectwise' definitions of the
tensor and cotensor, this makes $\Func(\cC,\cV)$ into a closed
$\cV$-module category (see Appendix A for terminology).

If $M\colon \cC \ra \cV$ is a functor and $X\in \cV$, one notes that
there is a canonical isomorphism $T(M\tens X)\iso (TM)\tens X$; this
follows from the explicit description of $\sigma \tens (\blank)$ given
above.  The map of functors $M\tens \F(M,N) \ra N$ therefore gives
rise to a map $TM \tens \F(M,N) \ra TN$, or a map $\eta_{M,N}\colon
\F(M,N)\ra \F(TM,TN)$ by adjointness.

If $M$ and $N$ are $\sigma$-modules then they come equipped with maps
of functors $TM \ra M$ and $TN \ra N$.  One defines $\F_\sigma(M,N)\in
\cV$ as the equalizer of the two obvious maps $\F(M,N) \dbra \F(TM,N)$
(to define one of the maps one uses $\eta_{M,N}$).  With this
definition---as well as the objectwise definitions for the tensor and
cotensor---the category of $\sigma$-modules becomes a closed
$\cV$-module category.  The adjunction $(\S_\cC-mod) \adjoint
(\sigma-mod)$ is a $\cV$-adjunction.  Using this together with the
observation that $\sigma(a,\blank)=\sigma \tens \S_\cC(a,\blank)$,
one sees that there are natural isomorphisms
$\F_\sigma(\sigma(a,\blank),M) \iso M(a)$.

\begin{prop}
\label{pr:module}
Assume $\cC$ is small and $\cV$ is a combinatorial, symmetric monoidal model
category satisfying the monoid axiom.  Let $\sigma$ be an enrichment
of $\cM$ by $\cV$.  Then there is a cofibrantly-generated model
structure on the category of left $\sigma$-modules in which a map
$M\ra M'$ is a weak equivalence or fibration precisely when $M(a) \ra
M'(a)$ is a weak equivalence or fibration for every $a\in \cC$.  This
makes the category of left $\sigma$-modules into a $\cV$-model
category.  If the unit $S\in \cV$ is cofibrant, then the free modules
$\sigma(a,\blank)$ are cofibrant.
\end{prop}

\begin{proof}
Take the generating cofibrations (resp. trivial cofibrations) to be
maps $\sigma(a,\blank)\tens A \ra \sigma(a,\blank)\tens B$ where $A\ra
B$ is a generating cofibration (resp. trivial cofibration) of $\cV$
and $a\in \cC$ is any object.  Checking that this gives rise to a
cofibrantly-generated model structure is a routine application of
\cite[Th. 11.3.1]{H}.  The other statements are routine verifications
as well.  See also \cite[Th. A.1.1]{SS2}.
\end{proof}

\begin{remark}
Of course everything above also works for {\it right\/} $\sigma$-modules.
\end{remark}

\subsection{Bimodules}
\label{se:bimodules2}
Suppose $\sigma$ is an enrichment of $\cC$ by $\cV$, and $\tau$ is an
enrichment of $\cD$ by $\cV$.
Define $\sigma \tens \tau$ to be the enrichment on $\cC \times \cD$
given by $(\sigma\tens \tau)((c_1,d_1),(c_2,d_2)) = \sigma(c_1,c_2)
\tens \tau(d_1,d_2)$.  Define $\sigma^{op}$ to be the enrichment of
$\cC^{op}$ given by $\sigma^{op}(a,b)=\sigma(b,a)$.  Finally, define a
$\sigma-\tau$ bimodule to be a left $\tau^{op}\tens\sigma$-module.

\begin{remark}
Upon unraveling the above definition, the reader will find that it is
equivalent with the more naive (and concrete) version given in
Section~\ref{se:modelenrich} for the case $\cC=\cD$.  The notational
conventions of that naive definition dictated the use of
$\tau^{op}\tens \sigma$ rather than $\sigma\tens \tau^{op}$ in the
above definition.
\end{remark}

It follows from Proposition~\ref{pr:module} that the category of
$\sigma-\tau$ bimodules has a model structure in which weak
equivalences and fibrations are determined objectwise.

Note that if $M$ is a $\sigma-\tau$ bimodule, then for any $a\in \cC$
the functor $M(a,\blank)$ is a left $\sigma$-module and the functor
$M(\blank,a)$ is a right $\tau$-module.

\subsection{The main proof} 
Exactly following \cite[Lem. A.2.3]{SS2}, we can now conclude the

\begin{proof}[Proof of Proposition~\ref{pr:equiv=direct}]
We will sketch the proof for the reader's convenience.  Suppose
$\sigma$ and $\tau$ are model enrichments of $\cM$ by $\cV$, defined
over some small category $\cI$ consisting of cofibrant-fibrant
objects.  Assume there is a quasi-equivalence between them given by
the pointed bimodule $M$.  If the composites $\sigma(a,b)\tens S \ra
\sigma(a,b) \tens M(a,a) \ra M(a,b)$ are all trivial fibrations (or if
the corresponding maps $\tau(a,b) \ra M(a,b)$ are all trivial
fibrations) then the proof is exactly as in [loc. cit].

For the general case, we first replace $M$ with a fibrant model in the
category of $\sigma-\tau$ bimodules over $\cI$; this makes $M$
objectwise fibrant.  For each $a\in \cI$, the distinguished map $S\ra
M(a,a)$ gives a map of right $\tau$-modules $\cF_a=\tau(\blank,a) \ra
M(\blank,a)$.  We apply our functorial factorization in the model
category of right $\tau$-modules to obtain $\cF_a \cofib N_a \trfib
M(\blank,a)$.  As the factorization is functorial, for every map $a\ra
b$ in $\cI$ there is an induced map of right $\tau$-modules $N_a \ra
N_b$.  Note that each $N_a$ is both cofibrant and fibrant as a
$\tau$-module: the fibrancy is immediate, but the cofibrancy uses that
$\cF_a$ is cofibrant (which in turn depends on the unit $S\in \cV$
being cofibrant).  Let $\cE$ be the model enrichment of $\cI$ given by
$\cE(a,b)=\F_\tau(N_a,N_b)$.

Define $U$ to be the $\sigma-\cE$ bimodule
$U(a,b)=\F_\tau(N_a,M(\blank,b))$ and $W$ to be the $\cE-\tau$
bimodule given by $W(a,b)=\F_\tau(\cF_a,N_b)$.  The fact that $W$ is a
right $\tau$-module uses the existence of maps $\tau(i,j) \ra
\F_\tau(\cF_i,\cF_j)$, which is easily established.  
One sees that $U$ and
$W$ are naturally pointed, and give quasi-equivalences between $\sigma$
and $\cE$, and between $\cE$ and $\tau$, respectively.  Moreover, we
are now in the case handled by the first paragraph of this proof,
because for $U$ and $W$ the appropriate maps are trivial fibrations.
So we get a zig-zag of four direct equivalences between $\sigma$ and
$\tau$.
\end{proof}


\section{The additive case: Homotopy enrichments over $\Spe^\Sigma(s\Ab)$}
\label{se:additive}

We'll say that a model category is \dfn{additive} if its underlying
category has a zero object and is enriched over abelian groups.  If
$\cM$ is an additive, stable, combinatorial model category, we will
produce a model enrichment of $\cM$ by $\Spe^\Sigma(s\Ab)$.  This
allows us in particular to attach to every object $X\in \cM$ an
isomorphism class in $\ho(\Ring[\Spe^\Sigma (s\Ab)])$.  Write
$\hdga(X)$ for any object in this isomorphism class.  

By \cite{S}, the homotopy category of $\Ring[\Spe^\Sigma(s\Ab)]$ is the
same as the homotopy category of dgas over $\Z$.  So $\hdga(X)$ can be
regarded as a `homotopy endomorphism dga'.  Unlike the homotopy
endomorphism spectra of Section~\ref{se:endo}, however, this dga is
{\it not\/} an invariant of Quillen equivalence.  It does act as an
invariant if one restricts to strings of Quillen equivalences
involving only
additive model categories, though.

\medskip

Here are the basic results (see (\ref{se:background}) for additional
terminology):

\begin{prop}
\label{pr:hdga=EM}
Given $X\in \cM$ as above, the homotopy endomorphism spectrum
$\hEnd(X)$ is the Eilenberg-MacLane spectrum associated to 
$\hdga(X)$.
\end{prop}

\begin{prop}
\label{pr:hdga-invariance}
Let $\cM$ and $\cN$ be additive, stable, combinatorial model
categories.  Suppose $\cM$ and $\cN$ are Quillen equivalent through a
zig-zag of additive (but not necessarily combinatorial) model
categories.  Let $X\in \cM$, and let $Y\in \ho(\cN)$ correspond to $X$
under the derived equivalence of homotopy categories.  Then $\hdga(X)$
and $\hdga(Y)$ are weakly equivalent in $\Ring[\Spe^\Sigma(s\Ab)]$.
\end{prop}

\begin{prop}
\label{pr:hdga-spectral}
Let $\cM$ be additive, stable, combinatorial, and an
$\Spe^\Sigma(s\Ab)$-model category.  Let $X\in \cM$ be
cofibrant-fibrant.  Then $\hdga(X)$ is weakly equivalent
to the cotensor object $\F(X,X)$.
\end{prop}

The proofs of the above two results are for the most part similar to
the corresponding results for homotopy endomorphism spectra.  One
difference is that they depend on developing a theory of universal
{\it additive\/} model categories.  Another, more important,
difference is the following.  Recall from Proposition~\ref{pr:univ}(b)
that any zig-zag of Quillen equivalences between two pointed model
categories (with the intermediate steps not necessarily pointed) could
be replaced by a simple zig-zag where the third model category is also
pointed.  In contrast to this, it is not generally true that a zig-zag
of Quillen equivalences between two {\it additive\/} model categories
(with intermediate steps not necessarily additive) can be replaced by
a simple zig-zag where the middle step is also additive.  This is only
true if we assume that all the intermediate steps are additive in the
first place.

\medskip

\subsection{Background}
\label{se:background}
If $\cM$ is a monoidal model category which is combinatorial and
satisfies the monoid axiom, then by \cite[Th. 4.1(3)]{SS1} the
category of monoids in $\cM$ has an induced model structure where the
weak equivalences and fibrations are the same as those in $\cM$.
We'll write $\Ring[\cM]$ for this model category.  If $\cN$ is another
such monoidal model category and $L\colon \cM \adjoint \cN \colon R$
is a Quillen pair which is weak monoidal in the sense of
\cite[Def. 3.6]{SS3}, then there is an induced Quillen map
$\Ring[\cM] \ra \Ring[\cN]$.  This is a Quillen equivalence if $\cM
\ra \cN$ was a Quillen equivalence and the units in $\cM$ and $\cN$
are cofibrant \cite[Th. 3.12]{SS3}.

The adjunction $\Set_* \adjoint \Ab$ is strong monoidal, and therefore
induces strong monoidal Quillen functors $\Spe^\Sigma(\sSet_*)\adjoint
\Spe^\Sigma(s\Ab)$.  Therefore one gets a Quillen pair $F\colon
\Ring[\Spe^\Sigma] \adjoint \Ring[\Spe^\Sigma(s\Ab)]\colon U$.  By the
\dfn{Eilenberg-MacLane ring spectrum} associated to an $R\in
\Ring[\Spe^\Sigma(s\Ab)]$ we simply mean the ring spectrum $UR$.

\subsection{Universal additive model categories}
Let $\cC$ be a small, semi-additive category.  This means the Hom-sets
of $\cC$ have a natural structure of abelian groups, and $\cC$ has a
zero-object \cite[VIII.2]{M}---the `semi' is to indicate that $\cC$
need not have direct sums.  One says that a functor $F\colon
\cC^{op}\ra \Ab$ is \dfn{additive} if $F(0)\iso 0$ and for any two
maps $f,g\colon X \ra Y$ in $\cC$ one has $F(f+g)=F(f)+F(g)$.
Note that for every $X\in \cC$, the representable functor $rX$ defined
by $U\mapsto \cC(U,X)$ is additive.  

Let $\Func(\cC^{op},\Ab)$ denote
the category of all functors.  The Yoneda Lemma does not hold in this
category: that is, if $F\in \Func(\cC^{op},\Ab)$ one need not have
$\Hom(rX,F)\iso F(X)$ for all $X\in \cC$.  But it is easy to
check that this {\it does\/}
hold when $F$ is an additive functor.

Let $\Func_{ad}(\cC^{op},\Ab)$ denote the full subcategory of additive
functors.  The following lemma records several basic facts about this
category, most of which follow from the Yoneda Lemma.

\begin{lemma} 
\label{le:functorcats}
Let $\cC$ be a small, semi-additive category.
\begin{enumerate}[(a)]
\item Colimits and limits in $\Func_{ad}(\cC^{op},\Ab)$ are the same
as those in $\Func(\cC^{op},\Ab)$.  
\item Every additive functor $F\in \Func(\cC^{op},\Ab)$ is isomorphic
to its canonical colimit with respect to the embedding $r\colon \cC
\inc \Func(\cC^{op},\Ab)$.  That is, the natural map
$\colim\limits_{rX \ra F} rX \ra F$ is an isomorphism.
\item
The additive functors in $\Func(\cC^{op},\Ab)$ are precisely those
functors which are colimits of representables.
\item The inclusion $i\colon \Func_{ad}(\cC^{op},\Ab) \inc
\Func(\cC^{op},\Ab)$ has a left adjoint $\Ad$ (for `additivization'),
and the composite $\Ad \circ i$ is naturally isomorphic to the
identity.
\item Suppose given a co-complete, additive category $\cA$ and an
additive functor $\gamma\colon \cC \ra \cA$.  Define $\Sing\colon \cA \ra
\Func_{ad}(\cC^{op},\Ab)$ by letting $\Sing(a)$ be the
functor $a\mapsto \cA(\gamma c,a)$.  Then $\Sing$ has a left adjoint
$\Real$, and there are natural isomorphisms $\Real(rX) \iso
\gamma(X)$.
\end{enumerate}
\end{lemma}

\begin{proof}
Left to the reader.  
\end{proof}

By \cite[Th. 11.6.1]{H} the category $\Func(\cC^{op},s\Ab)$ has a
cofibrantly-generated model structure in which the weak equivalences
and fibrations are defined objectwise.  We will need the analagous
result for the category of additive functors:

\begin{lemma}
\label{le:functor-mc}
The category $\Func_{ad}(\cC^{op},s\Ab)$ has a cofibrantly-generated
model structure in which the weak equivalences and fibrations are
defined objectwise.  This model structure is simplicial, left proper,
and cellular.
\end{lemma}

\begin{proof}
The proof uses the adjoint pair $(\Ad,i)$ to create the model
structure, as in \cite[Th. 11.3.2]{H}.  Recall that the model category
$\Func(\cC^{op},s\Ab)$ has generating trivial cofibrations $J=\{rX
\tens \Z[\Lambda^{n,k}] \ra rX \tens \Z[\Delta^n] \mid X\in \cC\}$.  Our
notation is that if $K\in \sSet$ then $\Z[K] \in s\Ab$ is the
levelwise free abelian group on $K$; and if $A\in s\Ab$ then $rX\tens
A$ denotes the presheaf $U \mapsto \cC(U,X) \tens A$ (with the tensor
performed levelwise).

To apply \cite[11.3.2]{H} we must verify that the functor $i$ takes
relative $\Ad(J)$-cell complexes to weak equivalences.  However, note
that the domains and codomains of maps in $J$ are all additive
functors (since representables are additive), and so $\Ad(J)=J$.  The
fact that forming pushouts in $\Func_{ad}(\cC^{op},s\Ab)$ and
$\Func(\cC^{op},s\Ab)$ give the same answers (by
Lemma~\ref{le:functorcats}(a)) therefore shows that the $\Ad(J)$-cell
complexes are indeed weak equivalences in $\Func(\cC^{op},s\Ab)$. 

Finally, it is routine to check that the resulting model structure is
simplicial, left proper, and cellular.
\end{proof}

From now on we will write $\Ua \cC$ for the category
$\Func_{ad}(\cC^{op},s\Ab)$ with the model structure provided by the
above lemma.  The reason for the notation is provided by the next result.

\begin{thm} 
\label{th:additive-univ}
Let $\cM$ be an additive model category.
\begin{enumerate}[(a)]
\item Suppose $\cC$ is a small, semi-additive category and $\gamma\colon
\cC \ra \cM$ is an additive functor.  Then there is a Quillen pair
$\Real \colon \Ua\cC \adjoint \cM \colon \Sing$ together with a
natural weak equivalence $\Real \circ r \we \gamma$.  
\item If $\cM$ is combinatorial then there
is a Quillen equivalence $\Ua\cC/S \we \cM$ for some small, semi-additive
category $\cC$ and some set of maps $S$ in $\Ua\cC$.
\item Suppose $\cM \bwe \cM_1 \we \cdots \bwe \cM_n \we \cN$ is a
zig-zag of Quillen equivalences in which all the model categories are
additive.  If $\cM$ is combinatorial, there is a simple zig-zag of
equivalences 
\[ \cM \bwe \Ua\cC/S \we \cN 
\]
such that the derived equivalence $\ho(\cM)\he \ho(\cN)$ is
isomorphic to the derived equivalence given by the original zig-zag.
\end{enumerate}
\end{thm}

\begin{proof}
The proofs for (a) and (c) are simple, and exactly follow the case for
$\U\cC$ (see \cite[Prop. 2.3, Cor. 6.5]{D1}).  The proof of (b) is
slightly more complicated, and will be postponed until the end of this
section.
\end{proof}

\begin{remark}
The result in (c) is false if one does not assume that all the
$\cM_i$'s are additive.  For an example, let $R$ be the dga $\Z[e;
de=2]/(e^4)$ and let $T$ be the dga $\Z/2[x;dx=0]/(x^2)$, where both
$e$ and $x$ have degree $1$.  Let $\cM$ and $\cN$ be the categories of
$R$- and $T$-modules, respectively.  These turn out to be Quillen
equivalent, but they cannot be linked by a zig-zag of Quillen
equivalences between {\it additive\/} model categories.  A
verification of these claims can be found in \cite[Example 6.10]{DS}.
\end{remark}

\subsection{Endomorphism objects}
Let $\cM$ be an additive, stable, combinatorial model category.  By
Theorem~\ref{th:additive-univ} there is a Quillen equivalence $\Ua
\cC/S \ra \cM$ for some small, additive category $\cC$ and some set of
maps $S$ in $\Ua\cC$.  The category $\Ua\cC/S$ is simplicial, left
proper, and cellular, so we may form $\Spe^\Sigma (\Ua\cC/S)$. 
Since $\Ua\cC/S$ is stable (since $\cM$ was), we obtain a zig-zag of
Quillen equivalences
\[ \cM \bwe \Ua\cC/S  \we \Spe^\Sigma(\Ua\cC/S).
\]

The category $\Ua\cC$ is an $s\Ab$-model category, and therefore
$\Spe^\Sigma(\Ua\cC/S)$ is an $\Spe^\Sigma(s\Ab)$-model
category.  We can transport this enrichment onto $\cM$ via the Quillen
equivalences, and therefore get an element $\sigma_\cM\in
\ME_0(\cM,\Spe^\Sigma(s\Ab))$.  Just as in Section~\ref{se:main}, one
shows that this quasi-equivalence class does not depend on the choice
of $\cC$, $S$, or the Quillen equivalence $\Ua\cC/S \we \cM$.

Let $X\in \cM$, and let $\tilde{X}$ be a cofibrant-fibrant object
weakly equivalent to $X$.  We write \mdfn{$\hdga(X)$} for any object
in $\Ring[\Spe^\Sigma(s\Ab)]$ having the homotopy type of
$\sigma_\cM(\tilde{X},\tilde{X})$, and we'll call this the \mdfn{additive
homotopy endomorphism object of $X$}.  By
Corollaries~\ref{co:endo=invariant} and \ref{co:endo=coffib}, this
homotopy type depends only on the homotopy type of $X$ and the
quasi-equivalence class of $\sigma_\cM$---and so it is a well-defined
invariant of $X$ and $\cM$.

\begin{proof}[Proof of Proposition~\ref{pr:hdga-invariance}]
This is entirely similar to the proof of Theorem~\ref{th:invariance}.
\end{proof}

\begin{proof}[Proof of Proposition~\ref{pr:hdga-spectral}]
Same as the proof of Proposition~\ref{pr:spectralenrich}.
\end{proof}

\begin{proof}[Proof of Proposition~\ref{pr:hdga=EM}]
We know that there exists a zig-zag of Quillen equivalences $\cM \bwe
\Ua\cC/S \we \Spe^\Sigma(\Ua\cC/S)$.  Therefore, using
Theorem~\ref{th:invariance} and Proposition~\ref{pr:hdga-invariance}
we may as well assume $\cM=\Spe^\Sigma(\Ua \cC/S)$.  This is an
$\Spe^\Sigma(s\Ab)$-model category, and so for any object $X$ we have a
ring object $F(X,X)$ in $\Spe^\Sigma(s\Ab)$.  The adjoint functors
$\Set_* \adjoint \Ab$ induce a strong monoidal adjunction $F\colon
\Spe^\Sigma(\sSet_*) \adjoint \Spe^\Sigma(s\Ab)\colon U$.  The
$\Spe^\Sigma(s\Ab)$-structure on $\cM$ therefore yields an induced
$\Spe^\Sigma$-structure as well (see Section~\ref{se:lifting}).  In
this structure, the endomorphism ring spectrum of $X$ is precisely
$U[F(X,X)]$.  Proposition~\ref{pr:spectralenrich} tells us this has
the homotopy type of the ring spectrum $\hEnd(X)$, at least when $X$
is cofibrant-fibrant.  And Proposition~\ref{pr:hdga-spectral} says
that $F(X,X)$ has the homotopy type of $\hdga(X)$.  This is all we
needed to check.
\end{proof}

\subsection{Additive presentations}
We turn to the proof of Theorem~\ref{th:additive-univ}(b).  This will
be deduced from the work of \cite{D2} plus some purely formal
considerations.

Let $\cM$ be a combinatorial model category.  By \cite[Prop. 3.3]{D2},
there is a small category $\cC$ and a functor $\cC \ra \cM$ such that
the induced map $L\colon \U\cC \ra \cM$ is {\it homotopically
surjective\/} (see \cite[Def. 3.1]{D2} for the definition).  Then
\cite[Prop. 3.2]{D2} shows that there is a set of
maps $S$ in $\U\cC$ which the derived functor of $L$ takes to weak
equivalences, and such that the resulting map $\U\cC/S \ra \cM$ is a
Quillen equivalence.

Now suppose that $\cM$ was also an additive model category.  By
examining the proof of \cite[Prop. 3.3]{D2} one sees that $\cC$ may be
chosen to be a semi-additive category and the functor $\gamma\colon
\cC \ra \cM$ an additive functor (the category $\cC$ is a certain full
subcategory of the cosimplicial objects over $\cM$).  By
Theorem~\ref{th:additive-univ}(a) there is an induced map
$F\colon\Ua\cC \ra \cM$.  Again using \cite[Prop. 3.2]{D2}, it will be
enough to prove that this map is homotopically surjective.

Consider now the following sequence of adjoint pairs:
\[ \xymatrix{
\Func(\cC^{op},\sSet) \ar@<0.5ex>[r]^{\Z}
& \Func(\cC^{op},s\Ab) \ar@<0.5ex>[r]^{\Ad}\ar@<0.5ex>[l]^U
& \Func_{ad}(\cC^{op},s\Ab) \ar@<0.5ex>[r]^-F\ar@<0.5ex>[l]^i
& \cM \ar@<0.5ex>[l]
}
\]
The composite of the right adjoints is clearly the right adjoint of
$L$, so the composite of the left adjoints is $L$.  We have
constructed things so that this composite is homotopically surjective,
and we are trying to show that $F$ is also homotopically surjective.

\begin{lemma}
\label{le:AdZ}
If $X\in \cC$ then $\Ad(\Z (rX))\iso rX$ 
(or to be more precise, $Ui(\Ad(\Z (rX))) \iso rX$).
\end{lemma}

\begin{proof}
This is clear, since the functors $F\mapsto \Func_{ad}(\Ad(\Z(rX)),F)$ and
$F\mapsto \Func_{ad}(rX,F)$ are both naturally isomorphic to $F(X)$.   
\end{proof}

Let $G \in \Func_{ad}(\cC^{op},s\Ab)$.  Let $QG$ be the simplicial
presheaf whose $n$th level is
\[ \coprod_{rX_n \ra rX_{n-1} \ra \cdots \ra rX_0 \ra G_n} (rX_n)
\]
where the coproduct is in $\Func(\cC^{op},\sSet)$.  The simplicial
presheaf $QG$ is treated in detail in \cite[Sec. 2.6]{D1}, as it is a
cofibrant-replacement functor for $\U\cC$.   Likewise,
let $Q_{ad}G$ be the simplicial presheaf whose $n$th level is
\[ \bigoplus_{rX_n \ra rX_{n-1} \ra \cdots \ra rX_0 \ra G_n} (rX_n)
\]
where the coproduct is now in $\Func(\cC^{op},s\Ab)$.  The proof of
\cite[Prop. 2.8]{D2} showing that $Q$ is a cofibrant-replacement
functor for $\U\cC$ adapts verbatim to show that $Q_{ad}$ is a
cofibrant-replacement functor for $\Ua\cC$.  Note that by
Lemma~\ref{le:AdZ} we have $Q_{ad}G=\Ad(\Z (QG))$, since
$\Ad$ and $\Z(\blank)$ are left adjoints and therefore
preserve coproducts.  

Finally we are in a position to conclude the

\begin{proof}[Proof of Theorem~\ref{th:additive-univ}(b)]
We have reduced to showing that $F\colon \Ua\cC \ra\cM$ is
homotopically surjective.  Let $\Sing$ be the right adjoint of $F$.
Then we must show that for every fibrant object $X\in \cM$ the induced
map $FQ_{ad}(\Sing X)\ra X$ is a weak equivalence.

However, the fact that $L\colon \U\cC \ra \cM$ is homotopically
surjective says that $LQ (Ui\Sing X) \ra X$ is a weak equivalence in
$\cM$.  And we have seen above that
\[ FQ_{ad}Ui\Sing X = F \bigl [\Ad (\Z (Q \Sing X))\bigr ] = L Q\Sing X, 
\]
so we are done.
\end{proof}


\appendix

\vspace{0.3in}

\section{$\cD$-model categories}

In the body of the paper we need to deal with {\it spectral\/} model
categories.  These are model categories which are enriched,
tensored, and cotensored over the model category of symmetric spectra,
and where the analog of SM7 holds.  In this appendix we briefly review
some very general material relevant to this situation.  We assume the
reader already has some experience in this area (for instance in the
setting of {\it simplicial\/} model categories), and for that reason
only give a broad outline.

\subsection{Basic definitions}
Let $\cD$ be a closed symmetric monoidal category.  The `symmetric
monoidal' part says we are given a bifunctor $\tens$, a unit object
$\unit_{\cD}$, together with associativity, commutativity, and unital
isomorphisms making certain diagrams commute (see \cite[Defs. 4.1.1,
4.1.4]{Ho2} for a nice summary).  The `closed' part says that there is
also a bifunctor $(d,e) \mapsto \mD(d,e)\in \cD$ together with a natural
isomorphism 
\[ \cD(a, \mD(d,e)) \iso \cD(a\tens d, e).
\]
Note that, in particular, this gives us isomorphisms
$\cD(\unit_{\cD},\mD(d,e)) \iso \cD(\unit_{\cD}\tens d,e) \iso \cD(d,e)$.

We define a \mdfn{closed $\cD$-module category} to be a category $\cM$
equipped with natural constructions which assign
to every $X,Z\in \cM$ and $d\in \cD$ objects
\[ X\tens d \in \cM, \qquad \F(d,Z) \in \cM, \qquad \text{and} \qquad 
\mM_{\cD}(X,Z) \in \cD.
\]
One requires, first, that there are natural isomorphisms $(X\tens
d)\tens e\iso X\tens (d\tens e)$ and $X\tens \unit_{\cD} \iso X$
making certain diagrams commute (see \cite[Def. 4.1.6]{Ho2}).  One
also requires natural isomorphisms
\begin{equation}
\label{eq:adjoint}
 \cM(X\tens d,Z) \iso \cM(X, \F(d,Z)) \iso \cD(d, \mM_{\cD}(X,Z))
\end{equation}
(see \cite[4.1.12]{Ho2}).

\begin{remark}
\label{re:iso}
Taking $d=\unit_{\cD}$, note that we obtain isomorphisms
$\cM(X,Z)\iso \cM(X\tens \unit,Z) \iso
\cD(\unit,\mM_{\cD}(X,Z))$.  
\end{remark}

\begin{prop}
Suppose $\cD$ is a symmetric monoidal category, and $\cM$ is a
closed $\cD$-module category.  Then one 
has canonical
isomorphisms
\[ \mM_{\cD}(X\tens d,Z) \iso \mM_{\cD}(X,\F(d,Z)) \iso \mD(d,
\mM_{\cD}(X,Z))
\]
of objects in $\cD$.  Applying $\cD(\unit_\cD,\blank)$ to these
isomorphisms yields the isomorphisms in (\ref{eq:adjoint}).
\end{prop}

\begin{proof}
The Yoneda Lemma says that two objects $a,b \in \cD$ are isomorphic if
and only if there is a natural isomorphism $\cD(e,a) \iso \cD(e,b)$,
for $e\in \cD$.  The proof of the proposition is straightforward using
this idea.
\end{proof}

\begin{prop}
Suppose $\cD$ is a symmetric monoidal category, and $\cM$ is a closed
$\cD$-module category.  Then there are `composition' maps
\[ \mM_{\cD}(Y,Z)\tens \mM_{\cD}(X,Y) \ra \mM_{\cD}(X,Z), \]
natural in $X$, $Y$, and $Z$.  These maps satisfy associativity and
unital conditions.  The induced map
\[ \mD(\unit,\mM_{\cD}(Y,Z)) \tens \mD(\unit,\mM_{\cD}(X,Y)) \ra
\mD(\unit,\mM_{\cD}(X,Z))
\]
coincides with the composition in $\cM$ under the isomorphisms from
Remark~\ref{re:iso}.
\end{prop}

\begin{proof}
We will only construct the maps, leaving the other verifications to
the reader.  The adjointness isomorphisms from (\ref{eq:adjoint})
give rise to natural maps $X\tens \mM(X,Y) \ra Y$ (adjoint to the
identity $\mM(X,Y)\ra \mM(X,Y)$).   There is a corresponding map
$Y\tens \mM(Y,Z) \ra Z$.  Now consider the composite
\[ X\tens [\mM(X,Y) \tens \mM(Y,Z)] \iso [X\tens \mM(X,Y)] \tens
\mM(Y,Z) \ra Y\tens \mM(Y,Z) \ra Z.
\]
Adjointness now gives $\mM(X,Y)\tens \mM(Y,Z) \ra \mM(X,Z)$, and
finally one uses that $\cD$ is symmetric monoidal.
\end{proof}

\begin{remark}
The basic definition of a $\cD$-module category doesn't really need
$\cD$ to be {\it symmetric\/} monoidal.  In fact, in \cite{Ho2} this is
not assumed.  However, the above propositions definitely need the
symmetric hypothesis.
\end{remark}

A \dfn{symmetric monoidal model category} consists of a closed
symmetric monoidal category $\cM$, together with a model structure on
$\cM$, satisfying two conditions:
\begin{enumerate}[(1)]
\item The analog of SM7, as given in either \cite[4.2.1]{Ho2} or
\cite[4.2.2(2)]{Ho2}. 
\item A unit condition given in \cite[4.2.6(2)]{Ho2}.
\end{enumerate}

Finally, let $\cD$ be a symmetric monoidal model category.  A
\mdfn{$\cD$-model category} is a model category $\cM$ which is also a
closed $\cD$-module category and where the two conditions from
\cite[4.2.18]{Ho2} hold: these are again the analog of SM7 and a unit
condition.  

\subsection{Lifting module structures}
\label{se:lifting}
Suppose that $\cC$ and $\cD$ are symmetric monoidal model categories,
and that $L\colon \cC \adjoint \cD \colon R$ is a Quillen pair.  One says
this adjunction is \dfn{strong symmetric monoidal} if there are
isomorphisms $L(\unit_\cC)\iso \unit_\cD$ and $L(X\tens Y)\iso LX
\tens LY$ compatible with the associativity, commutativity, and unital
isomorphisms in $\cC$ and $\cD$.  

\begin{lemma}
\label{le:lift}
Assume that $L\colon \cC \adjoint \cD \colon R$ is a strong symmetric
monoidal Quillen adjunction.  Let $\cM$ be a $\cD$-model category.
Then $\cM$ also becomes a $\cC$-model category by setting
\[ X\tens c=X\tens L(c), \quad \F_\cC(c,Y)=\F(Lc,Y), \quad
\text{and}\quad
\und{\cM}_\cC(X,Y)=R \bigl [\und{\cM}_{\cD}(X,Y) \bigr].
\]
\end{lemma}

\begin{proof}
Routine.
\end{proof}

\subsection{Spectral model categories}
\label{se:spectral}
Let $\Spe^\Sigma=\Spe^\Sigma(\sSet_+)$ be the usual category of
symmetric spectra \cite{HSS}.  This is a symmetric monoidal model
category.  We will call an $\Spe^\Sigma$-model category simply a
\dfn{spectral model category}.

Note that there are adjoint functors $\sSet_+ \adjoint \Spe^\Sigma$
where the left adjoint is $K\mapsto \Sigma^\infty(K)$ and the right
adjoint is $\Ev_0$, the functor sending a spectrum to the space in its
$0$th level.  The functor $\Sigma^\infty$ is called $F_0$ in
\cite{HSS}.  These functors are strong symmetric monoidal (see
\cite[2.2.6]{HSS}).  Therefore any spectral model category becomes an
$\sSet_+$-model category in a natural way, via Lemma~\ref{le:lift}.

The adjoint functors $\sSet \adjoint \sSet_+$ (which are also strong
monoidal) in turn show that any $\sSet_+$-model structure gives rise
to an underlying simplicial model structure.

\subsection{Diagram categories}
Let $I$ be a small category.  If $\cD$ is cofibrantly-generated, then
$\cD^I$ has a model structure in which the weak equivalences and
fibrations are defined objectwise.  If $X\in \cD^I$ and $d\in \cD$,
define the two objects
$X\tens d$, $\F(d,X) \in \cD^I$ as follows:
\[ X\tens d\colon i \mapsto X(i)\tens d, \qquad \F(d,X)\colon i\mapsto
\F(d,X(i)).
\]
Also, if $X,Z\in \cD^I$ define $\und{\cD^I}_{\cD}(X,Z) \in \cD$ to be
the equalizer of
\[ \prod_i \mD(X(i),Z(i)) \dbra \prod_{j\ra k} \mD(X(j),Z(k)).
\]

\begin{lemma}
Assume $\cD$ is a cofibrantly-generated, symmetric monoidal model
category.  With the above definitions, $\cD^I$ is a $\cD$-model
category.
\end{lemma}

\begin{proof}
Straightforward.
\end{proof}

\subsection{Adjunctions}
\label{se:Afunctor}

\begin{lemma}
\label{le:Dadjunct}
Let $\cM$ and $\cN$ be closed $\cD$-module categories, and let
$L\colon \cM \adjoint \cN \colon R$ be adjoint functors.  The
following are equivalent:
\begin{enumerate}[(a)]
\item There are natural isomorphisms $\mN_{\cD}(LX,Y) \iso
\mM_{\cD}(X,RY)$ which after applying $\cD(\unit_{\cD},\blank)$ reduce
to the adjunction $\cN(LX,Y)\iso \cM(X,RY)$.
\item There are natural isomorphisms $L(X\tens d) \iso L(X) \tens d$
which reduce to the canonical isomorphism for $d=\unit_{\cD}$.
\item There are natural isomorphisms $R(\F(d,Z))\iso \F(d,RZ)$ which
reduce to the canonical isomorphism when $d=\unit_{\cD}$.  
\end{enumerate}
\end{lemma}

\begin{proof}
Left to the reader.
\end{proof}

In the situation of the above lemma, we'll say that the adjoint pair
is a \mdfn{$\cD$-adjunction} between $\cM$ and $\cN$.  When $\cM$ and
$\cN$ are $\cD$-model categories we'll say that $\cM \ra \cN$ is a
\mdfn{$\cD$-Quillen map} (resp.  \mdfn{$\cD$-Quillen equivalence}) if
it is both a Quillen map (resp. Quillen equivalence) and a
$\cD$-adjunction.  In this paper we mostly need simplicial and
spectral Quillen functors, i.e. the cases where $D=\sSet$ or
$D=\Spe^\Sigma$.

\begin{remark}
\label{re:adjunct}
Note that in the situation of a $\cD$-adjunction one may form the
following composite, for any $A,B\in \cN$:
\[ \mN(A,B) \ra \mN(LRA,B) \llra{\iso} \mM(RA,RB). \]
Similarly, one has a natural map $\mM(X,Y) \ra \mN(LX,LY)$ for $X,Y\in
\cM$.  It is a routine exercise to check that the adjunction
isomorphism $\mN(LA,X)\llra{\iso} \mM(A,RX)$ is equal to the composite
$\mN(LA,X) \ra \mN(RLA,RX) \ra \mN(A,RX)$,
just as for ordinary adjunctions.
\end{remark}

\medskip

Let $\cD$ be a cofibrantly-generated, symmetric monoidal model
category, and let $\cM$ be a $\cD$-model category.  Suppose $I$ is a
small category and $\gamma \colon I \ra \cM$ is a functor.  Define
$\Sing \colon \cM \ra \Func(I^{op},\cD)$ by sending $X\in \cM$ to the
functor $i \mapsto \mM_{\cD}(\gamma(i),X)$.  This has a left adjoint
$\Real\colon \Func(I^{op},\cD) \ra \cM$ which sends a functor $A$ to the
coequalizer
\[ \coprod_{j\ra k} \gamma(j) \tens A(k) \dbra \coprod_{i}
\gamma(i)\tens A(i).
\]

\begin{prop}
\label{pr:diagram-adjunct}
The adjoint pair $\Real\colon \Func(I^{op},\cD) \adjoint \cM \colon \Sing$
is a $\cD$-adjunction.
\end{prop}

\begin{proof}
One readily checks condition (c) in Lemma~\ref{le:Dadjunct}.
\end{proof}


\vspace{0.2in}

\section{Stabilization and localization}

Let $\cM$ be an $\sSet_+$-model category which is pointed, left
proper, and cellular.  Under these conditions
one may form the stabilized model category $\Spe^\Sigma \cM$
\cite{Ho1}, and this is again a left proper and cellular model
category.  Recall that there are Quillen pairs $F_i\colon \cM \adjoint
\Spe^\Sigma \cM\colon \Ev_i$, for every $i\geq 0$ ($F_0 X$ is also
written $\Sigma^\infty X$, and $F_iX$ is morally the $i$th
desuspension of $F_0 X$).

If $S$ is a set of maps between cofibrant objects in $\cM$, let
\[ 
S_{stab}=\{ F_i(A) \ra F_i(B) \, | \, A\ra B \in S \ \text{and}\  i\geq 0\}.
\]
Our goal is the following basic result about commuting stabilization
and localization:

\begin{prop}
\label{pr:locstab}
In the above situation, the model categories $\Spe^\Sigma (\cM/S)$ and
$(\Spe^\Sigma \cM)/S_{stab}$ are identical.
\end{prop}

\begin{proof}
The stable model structure on $\Spe^\Sigma \cM$ is formed in two steps.
One starts with the projective model structure $\Spe^\Sigma_{proj}\cM$
where fibrations and weak equivalences are levelwise (and cofibrations
are forced).  Then one localizes this projective structure at a
specific set of maps given in \cite[Def. 8.7]{Ho1}.  Call this set
$T_\cM$.  It is important that $T_\cM$ depends only on the generating
cofibrations of $\cM$.

So $\Spe^\Sigma(\cM/S)$ is the localization of $\Spe^\Sigma_{proj}
(\cM/S)$ at the set $T_{\cM/S}$.  Likewise, $(\Spe^\Sigma
\cM)/S_{stab}$ is the localization of $(\Spe^\Sigma_{proj}
\cM)/S_{stab}$ at the set of maps $T_{\cM}$.  But as the generating
cofibrations of $\cM$ and $\cM/S$ are the same, we have
$T_{\cM}=T_{\cM/S}$.  In this way we have reduced the proposition to
the statement that the model structures $\Spe^\Sigma_{proj} (\cM/S)$
and $(\Spe^\Sigma_{proj} \cM)/S_{stab}$ are identical.

The trivial fibrations in a model category and its Bousfield
localization are always the same.  This shows that the trivial
fibrations in the following categories are the same:
\[ \Spe^\Sigma_{proj}(\cM/S), \qquad
\Spe^\Sigma_{proj} \cM,\qquad
(\Spe^\Sigma_{proj}\cM)/S_{stab}.
\]
An immediate corollary is that the cofibrations are also the same in
these three model categories.  Note also that these are all simplicial
model categories, with simplicial structure induced by that on
$\cM$---and in particular that the simplicial structures are
identical.

Since the trivial fibrations in $\Spe^\Sigma_{proj}(\cM/S)$ and
$(\Spe^\Sigma_{proj}\cM)/S_{stab}$ are the same, it will suffice to
show that trivial cofibrations are also the same.  But a cofibration
$A\cof B$ is trivial precisely when the induced map on simplicial
mapping spaces $\Map(B,X) \ra \Map(A,X)$ is a weak equivalence for
every fibrant object $X$.  Since the model categories have the same
simplicial structures, we have reduced to showing that they have the
same class of fibrant objects.

A fibrant object in $\Spe^\Sigma_{proj}(\cM/S)$ is a spectrum $E$ such
that each $E_i$ is fibrant in $\cM/S$; this means $E_i$ is fibrant in
$\cM$, and for every $A\ra B$ in $S$ the induced map $\Map(B,E_i) \ra
\Map(A,E_i)$ is a weak equivalence (recall that $S$ consists of maps
between cofibrant objects).

A fibrant object in $(\Spe^\Sigma_{proj}\cM)/S_{stab}$ is a fibrant
spectrum $E\in \Spe^\Sigma_{proj} \cM$ (meaning only that each $E_i$ is
fibrant in $\cM$) which is $S_{stab}$-local.  The latter condition
means that for every $A\ra B$ in $S$ and for every $i$, the map
$\Map(F_i(B),E) \ra \Map(F_i(A),E)$ is a weak equivalence.  But the
adjoint pair $(F_i,\Ev_i)$ is a simplicial adjunction---one readily
checks condition (b) or (c) of Lemma~\ref{le:Dadjunct}.  So we have
$\Map(F_i(B),E)\iso \Map(B,\Ev_i(E))$, and the same for $A$.  This
verifies that the two classes of fibrant objects are the same, and
completes the proof.
\end{proof}


\bibliographystyle{amsalpha}

\end{document}